\input ./style/arxiv-general.cfg
\documentclass[aap,MSNbibl,dvips]{arximspdf}
\makeatletter
   \@ifpackageloaded{graphicx}{}{\usepackage{graphicx}}
\makeatother

\doi{10.1214/15-AAP1134}
\volume{26}
\issue{3}
\pubyear{2016}
\firstpage{1807}
\lastpage{1836}
\docsubty{FLA}

\makeatletter
\def\sfrac#1#2{#1/#2}

\newcommand{\wb}{\overline}
\newcommand{\wh}{\widehat}
\newcommand{\rrvert}{\vert}
\newcommand{\llvert}{\vert}
\renewcommand{\mid}{|}
\newcommand{\dotsc}{\dots}
\newcommand{\Q}{\mathbb Q}
\newcommand{\N}{\mathbb N}
\newcommand{\R}{\mathbb R}
\renewcommand{\P}{\mathbb P}
\newcommand{\E}{\mathbb E}
\newcommand{\bu}{\mathbf{u}}
\newcommand{\bv}{\mathbf{v}}
\newcommand{\prob}{\stackrel{P}{\rightarrow}}
\newcommand{\weak}{\Rightarrow}
\newcommand{\disteq}{\stackrel{d}{=}}
\newtheorem{theorem}{Theorem}[section]
\newtheorem{lemma}[theorem]{Lemma}
\newtheorem{prop}[theorem]{Proposition}
\newtheorem{cor}{Corollary}[theorem]
\newproclaim{remark}[theorem]{Remark}
\makeatother

\begin{document}
\begin{frontmatter}

\title{Hack's law in a drainage network model: A~Brownian web approach}
\runtitle{Hack's law in a drainage network model}

\begin{aug}
\author[A]{\fnms{Rahul}~\snm{Roy}\ead[label=e1]{rahul@isid.ac.in}},
\author[A]{\fnms{Kumarjit}~\snm{Saha}\corref{}\ead[label=e2]{kumarjitsaha@gmail.com}}
\and
\author[A]{\fnms{Anish}~\snm{Sarkar}\ead[label=e3]{anish@isid.ac.in}}
\runauthor{R. Roy, K. Saha and A. Sarkar}
\affiliation{Indian Statistical Institute}
\address[A]{Theoretical Statistics and Mathematics Unit\\
Indian Statistical Institute\\
7 S. J. S. Sansanwal Marg\\
New Delhi 110016\\
India\\
\printead{e1}\\
\phantom{E-mail: }\printead*{e2}\\
\phantom{E-mail: }\printead*{e3}}
\end{aug}

%
\received{\smonth{2} \syear{2015}}
%
\revised{\smonth{7} \syear{2015}}

%
\begin{abstract}
Hack [Studies of longitudinal stream profiles in Virginia and Maryland (1957). Report], while studying the
drainage system in the Shenandoah valley and the adjacent mountains of Virginia,
observed a power law relation $l \sim a^{0.6}$
between the length $l$ of a stream from its
source to a divide and the area $a$ of the basin that collects the
precipitation contributing to the stream as tributaries. We study the
tributary structure of Howard's drainage network model of headward
growth and branching studied
by Gangopadhyay, Roy and Sarkar [\textit{Ann. Appl. Probab.} \textbf{14} (2004) 1242--1266].
We show that the exponent of Hack's law is $2/3$ for Howard's model.
Our study is based on a scaling of the process whereby the limit of the
watershed area of a stream
is area of a Brownian excursion process.
To obtain this, we define a dual of the model and show that under
diffusive scaling,
both the original network and its dual converge jointly to the standard
Brownian web and its dual.
\end{abstract}

%
\end{frontmatter}

\section{Introduction}\label{Secintro}
River basin geomorphology is a very old subject of study initiated
by Horton \cite{H45}. 
Hack \cite{H57}, studying the
drainage system in the Shenandoah valley and the adjacent mountains of Virginia,
observed a power law relation
%
\begin{equation}
\label{Hacklaw} l \sim a^{0.6}
\end{equation}
between the length $l$ of a stream from its
source to a divide and the area of the basin $a$ that collects the
precipitation contributing to the stream as tributaries. 
Hack also
corroborated this power law
through the data gathered by Langbein \cite{L47} of nearly $400$ different streams
in the northeastern United States. This empirical relation~(\ref
{Hacklaw}) is
widely accepted nowadays
albeit with a different exponent (see Gray \cite{G61},
Muller \cite{M73})
and is called Hack's law. Mandelbrot \cite{M83} mentions Hack's law to
strengthen his contention
that ``\textit{if all rivers as well as their basins are mutually similar},\vspace*{1pt}
\textit{the fractal length-area argument
predicts} (\textit{river}'\textit{s length})$^{1/D}$ \textit{is proportional to}
(\textit{basin}'\textit{s area})$^{1/2} $'' where
$ D > 1 $ is the fractal dimension of the river.
In this connection, it is worth remarking that the Hurst exponent in fractional
Brownian motion and in time series analysis arose from the study of the Nile
basin by Hurst \cite{H27} where he proposed the relation $l_{\perp} =
l_{\parallel}^{0.9}$
as that governing the width, $l_{\perp}$, and the length,
$l_{\parallel}$, of
the smallest rectangular region containing the drainage system.

Various statistical models of drainage networks have been proposed
(see Rodriguez-Iturbe and Rinaldo \cite{RR97} for a detailed survey). In this paper, we study the
tributary structure of
a two-dimensional drainage network called the Howard's model of
headward growth
and branching (see Rodriguez-Iturbe and Rinaldo \cite{RR97}). Our study is based on a scaling of
the process and
we obtain the watershed area of a stream as the area of a Brownian
excursion process. This gives a statistical
explanation of Hack's law and justifies the remark of Giacometti et~al.
\cite{GMRR96}:
``\textit{From the results}, \textit{we suggest that a statistical framework referring
to the scaling invariance of the entire basin structure should be used in
the interpretation of Hack}'\textit{s law}.''

We first present an informal description of the model: suppose that the
vertices of
the $d$-dimensional lattice $ {\mathbb Z}^d$ are open or closed with
probability $ p$ $ (0 < p < 1) $
and $ 1-p$, respectively, independently of all other vertices.
Each open vertex $ \bu\in{\mathbb Z}^d$ represents a water source and
connects to
a unique open vertex $ \bv\in{\mathbb Z}^d$. These edges represent
the channels through which water can flow. 
The connecting vertex $ \bv$ is chosen so that the $d$th coordinate of
$ \bv$ is
one more than that of $ \bu$ and $ \bv$ has the minimum
$L_1$ distance from $ \bu$. In case of nonuniqueness of such a
vertex, we choose one of the closest
open vertices with equal probability, independently of everything else.

Let $V $ 
denote the set of \textit{open} vertices and $ h(\bu) $ denote the
uniquely chosen vertex to
which $ \bu$ connects, as described above. Set $ \langle\bu, h (\bu
) \rangle$ as the edge (channel) connecting
$ \bu$ and $ h(\bu) $.
From the construction, it follows that the random graph,
$ {\mathcal G} = (V, E) $ with edge set $ E:= \{ \langle\bu, h (\bu
) \rangle:
\bu\in V \}$, does not contain any circuit.
This model has been studied by
Gangopadhyay, Roy and Sarkar
\cite{GRS04} and
the following results were obtained.
%
\begin{theorem}
\label{GRS}
Let $0<p<1$.
\begin{longlist}[(ii)]
\item[(i)] For $d=2$ and $d=3$, ${\mathcal G}$ consists of one
single tree
almost surely,
and for $d\geq4$, ${\mathcal G}$ is a forest consisting of infinitely many
disjoint trees almost surely.
\item[(ii)] For any $d\geq2$, the graph ${\mathcal G}$ contains no
bi-infinite
path almost
surely.
\end{longlist}
\end{theorem}

In this paper, we consider only $d=2$. Before proceeding further, we
present a formal
description for $ d=2$ which will be used later. Fix $0< p < 1$ and let
$\{B_{\bu}: \bu= (\bu(1), \bu(2)) \in{\mathbb Z}^2 \}$ be
an i.i.d. collection of Bernoulli random variables with success
probability $p$. Set $ V = \{ \bu\in{\mathbb Z}^2: B_{\bu} = 1 \} $.
Let $\{U_{\bu}: \bu\in{\mathbb Z}^2 \}$ be
another i.i.d. collection of random variables, independent of the
collection of random variables $ \{ B_{\bu}: \bu\in{\mathbb Z}^2 \}$,
taking values in the set $ \{ 1, -1 \} $, with $ \P( U_{\bu} = 1 ) =
\P( U_{\bu} = -1) = 1/2 $.
For a vertex $ (x,t) \in{\mathbb Z}^2 $, we consider $ k_0 =
\min \{ \llvert  k \rrvert: k \in{\mathbb Z}, B_{(x + k, t+1)} = 1  \}
$. Clearly, $ k_0 $ is almost surely finite.
Now, we define
\begin{eqnarray*}
h (x,t):= \cases{ (x+k_0, t+1)\in V, &\quad if $(x -
k_0, t+1) \notin V$,
\cr
(x-k_0, t+1)\in V, &\quad if
$(x+k_0, t+1) \notin V$,
\cr
(x+U_{(x,t)}k_0, t+1)
\in V, &\quad if $(x\pm k_0, t+1) \in V$.}
\end{eqnarray*}

For any $ k \geq0 $, let
\begin{eqnarray*}
h^{k+1} (x, t) &:=& h \bigl( h^{k} (x, t) \bigr)\qquad
\mbox{with } h^0(x,t):= (x,t),
\\
C_k (x,t) &:=& \cases{ \bigl\{ (y, t-k) \in V: h^{k} ( y,
t-k) = (x,t) \bigr\}, &\quad if $(x,t) \in V$,
\cr
\varnothing, &\quad otherwise,}
\\
C(x,t) &:=& \bigcup_{ k \geq0 } C_k (x,t).
\end{eqnarray*}

Here,
$ h^k (x,t) $ represents the ``$k$th generation progeny'' of $ (x,t) $,
the sets $ C_k (x,t) $ and $ C (x,t) $ denote, respectively, the set of $
k$th generation ancestors and the set of all ancestors of $ (x,t) $; $
C(x,t) = \varnothing$ if
$ (x,t) \notin V$. In the
terminology of drainage network, $C (x,t)$ represents the region of
precipitation, the water from which is channelled through the open
point $ (x,t) $ (see Figure~\ref{figCx,t}). From Theorem \ref
{GRS}(ii), we have that $ C (x,t) $ is finite almost surely.

%
\begin{figure}[t]

\includegraphics{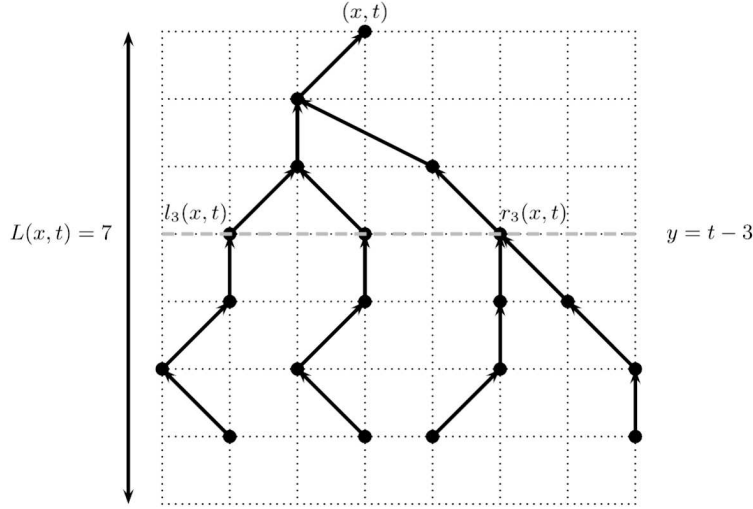}

\caption{The bold vertices on the line $y = t-3$ constitute the set
$C_3(x,t)$ and all the bold vertices together constitute the cluster $C(x,t)$.}\label{figCx,t}
\end{figure}


Now, we define
\[
L(x,t):= \inf\bigl\{ k \geq0: C_k (x,t) = \varnothing\bigr\},
\]
as the ``length of the channel,'' which as earlier is finite almost surely.
We observe that for any $(x,t)\in{\mathbb Z}\times{\mathbb Z}$,
$L(x,t) \geq0$ and the distribution
of $L(x,t)$ does not depend upon $(x,t)$. Our first
result is about the length of the channel. We remark here that
Newman, Ravishankar and Sun  \cite{NRS05} has a similar result in a set-up which allows crossing of paths.
%
\begin{theorem}
\label{clusterheight}
We have
\[
\lim_{ n \to\infty} \sqrt{ n } \mathbb{P} \bigl( L(0,0) > n \bigr) =
\frac{
1 }{
\gamma_0 \sqrt{\pi} },
\]
where $ \gamma_0^2:= \gamma_0^2 (p) = \frac{ (1-p)(2 - 2p +p^2) }{
p^2 ( 2-p)^2} $.
\end{theorem}

Next, we define
\begin{eqnarray*}
r_k (x,t) &:=& \cases{ \max\bigl\{ u: (u,t-k) \in C_k
(x,t) \bigr\}, &\quad if $0 \leq k < L(x,t)$, $(x,t) \in V$,
\cr
0, &\quad
otherwise,}
\\
l_k (x,t) &:=& \cases{ \min\bigl\{ u: (u,t-k) \in C_k
(x,t) \bigr\}, &\quad if $0 \leq k < L(x,t)$, $(x,t) \in V$,
\cr
0, &\quad
otherwise,}
\\
D_k(x,t) &:=& r_k (x,t) - l_k (x,t).
\end{eqnarray*}
The quantity $ D_k (x,t)$ denotes the \textit{width}
of the set of all $k$th generation ancestors of $ (x,t) $.
We define the \textit{width process} $
D_n^{(x,t)}(s) $ and the \textit{cluster process} $K_n^{(x,t)}(s)$
for $ s
\geq0$ as
follows: for $ k = 0, 1, \dots$ and $ k/ n
\leq s \leq
(k+1)/n $,
%
\begin{eqnarray}
\label{eqWidthClusterProcess} D_n^{(x,t)} (s) &: =& \frac{ D_k (x,t) }{ \gamma_0 \sqrt{n} } +
\frac{ (ns - [ns]) }{ \gamma_0 \sqrt{n} } \bigl( D_{k+1} (x,t) - D_k (x,t) \bigr),
\nonumber\\[-8pt]\\[-8pt]\nonumber
K_n^{(x,t)} (s) &: =& \frac{\# C_k (x,t) }{ \gamma_0 \sqrt{n} } + \frac{ (ns - [ns]) }{ \gamma_0 \sqrt{n} }
\bigl( \#C_{k+1} (x,t) - \#C_k (x,t) \bigr),
\end{eqnarray}
where $ \gamma_0 > 0 $ is as in the statement of Theorem \ref
{clusterheight}. In other words, $ D_n^{(x,t)} (s) $ [resp., $
K_n^{(x,t)} (s)$]
is defined $ D_k (x,t) / (\gamma_0 \sqrt{n} ) $ [resp., $\# C_k (x,t)
/ (\gamma_0 \sqrt{n})$] at time points $ s =
k/n $ and, at other time points defined by linear interpolation. %
The distributions of both $D_n^{(x,t)}$ and $K_n^{(x,t)}$ are
independent of
$(x,t)$.

To describe our results, we need to introduce two processes, Brownian meander
and Brownian excursion, studied by
Durrett, Iglehart and Miller
\cite{DIM77}. Let $ \{ W(s): s \geq0 \} $ be a standard Brownian
motion with $ W(0) = 0$. Let
$\tau_1:= \sup\{s \leq1: W(s)=0\}$ and $\tau_2:= \inf\{s \geq1:
W(s)=0\}$. Note that $\tau_1 < 1$ and $\tau_2 > 1$ almost surely.
The standard Brownian meander, $W^{+}(s)$, and the standard Brownian
excursion, $W^{+}_{0}(s)$, are given by
%
\begin{eqnarray}
\label{eqBMBE} W^{+}(s) &:=& \frac{\llvert   W(\tau_1 + s(1-\tau_1))\rrvert   }{\sqrt{1 - \tau
_1}}, \qquad s \in[0,1],
\\
W^{+}_{0}(s) &:=& \frac{\llvert  W(\tau_1 + s(\tau_2-\tau_1))\rrvert   }{\sqrt{\tau
_2-\tau_1}}, \qquad s\in[0,1].
\end{eqnarray}
Both of these processes are a continuous nonhomogeneous Markov process
(see Durrett and Iglehart \cite{DI77} and references therein).
Further, $W^{+}(0) = 0$ and, for $x \geq0$, $\P(W^{+}(1) \leq x) = 1
- \exp(-x^2/2)$, that is,
$W^{+}(1)$ follows a Rayleigh distribution.

We also need some random variables obtained as functionals of these two
processes. In particular, let
\begin{eqnarray*}
I^{+}_{0}&:=& \int_{0}^{1}W^{+}_{0}(t)\,dt
\quad\mbox{and}\quad M^{+}_{0}:= \max\bigl
\{W^{+}_{0}(t): t\in[0,1]\bigr\}. 
\end{eqnarray*}
%
Louchard and Janson  \cite{JL07} showed that, as $x \to\infty$, the distribution
function and the density are, respectively, given by
\begin{eqnarray*}
\P\bigl(I^{+}_0 &>& x
\bigr) \sim\frac{6 \sqrt{6}}{\sqrt{\pi}} x\exp{\bigl(-6x^2\bigr)} \quad\mbox{and}
\quad f_{I^+_0}(x) \sim\frac{72 \sqrt{6}}{\sqrt
{\pi}} x^2\exp{
\bigl(-6x^2\bigr)}.
\end{eqnarray*}
%
The random variable $M^{+}_0$ is continuous, having a strictly positive
density on $ (0, \infty)$ (see Durrett and Iglehart \cite{DI77}) and for $x > 0$,
\begin{eqnarray*}
\hspace*{-3pt}&&\P\bigl(M^{+}_0 \leq x
\bigr) = 1 + 2\sum^{\infty}_{k=1} \exp {
\bigl(-(2kx)^2/2\bigr)}\bigl[1 - (2kx)^2\bigr]\qquad
\mbox{with }\E\bigl(M^{+}_0\bigr) = \sqrt {\pi/2}.
\end{eqnarray*}
%

For $f \in C[0,\infty)$, let $f\mid _{[0,1]}$ denotes the restriction of
$f$ over $[0,1]$.
Our next result is about the weak convergence of the width process
$D^{(0,0)}_n\mid_{[0,1]}$ and
the cluster process $K^{(0,0)}_n\mid_{[0,1]}$ under diffusive scaling.
Here and subsequently, as is commonly used in statistics, we use the
notation $X\mid Y$ to denote the conditional random variable $X$
given $Y$.

%
\begin{theorem}
\label{BM}
As $n\rightarrow\infty$, we have:
\begin{longlist}[(ii)]
\item[(i)] $ D^{(0,0)}_n\mid _{[0,1]} \mid {\mathbf{1}}_{\{ L(0,0) > n
\}} \weak \sqrt{2} W^{+} $,
\item[(ii)] $\sup\{\llvert   pD_n^{(0,0)}(s) - K_n^{(0,0)}(s)\rrvert: s \in[0,1]\}
\mid {\mathbf{1}}_{\{ L(0,0) > n \}} \prob 0$.
\end{longlist}
\end{theorem}
%

The following corollary is an immediate consequence of Theorem \ref{BM}.
%
\begin{cor}
For $ u > 0 $, as $n \rightarrow\infty$ we have:
\begin{longlist}[(ii)]
\item[(i)]
$ \sqrt{ n} \mathbb{P} ( \# C_n (0,0) > \sqrt{n} \gamma_0 u )
\to \frac{ 1 }{ \gamma_0 \sqrt{\pi} } \exp( - u^2/4 p^2 )$,
\item[(ii)] $\mathbb{P} ( \sum_{k=0}^{n}\# C_k (0,0) > n^{3/2} \gamma_0 u
\mid  L(0,0) > n) \to \mathbb{P}(p\sqrt{2}I^{+} > u)$.
\end{longlist}
\end{cor}

Before we proceed to state Theorem \ref{BEarea}, we recall some
results regarding random vectors whose distribution functions have
regularly varying tails (see Resnick~\cite{R07}, page~172). A random vector $Z$
on $(0, \infty)^d$ with a distribution function $F$ has a regularly
varying tail if, as $n \to\infty$, there exists a sequence $b_n \to
\infty$ such that $n \P\{Z/b_n \in\cdot\} \stackrel{v}{\to} \nu
(\cdot)$ for some $\nu\in M_{+}$ where $M_{+}:= \{\mu: \mu$
is a nonnegative Radon measure on $(0, \infty)^d\}$. Here, $\stackrel
{v}{\to} $ denotes vague convergence.
It is in this context that Theorem \ref{BEarea} obtains a regularly
varying tail for the distribution of $(L(x,t),
(\#C(x,t))^{2/3})$; which justifies that the exponent of Hack's law is
$2/3$ for
Howard's model. In addition, we obtain a scaling law, with a Hack
exponent of $1/2$, for the length of the stream, vis-\`a-vis the maximum
width of the region of
precipitation, that is,
%
\begin{equation}
\label{eqndefMaxWidth} D_{\max}(0,0):= \max\bigl\{D_{k}(0,0):0\leq k
< L(0,0)\bigr\}.
\end{equation}
It should be noted that Leopold and Langbein
\cite{L62} obtained an exponent of $0.64$ through
computer simulations.
%
\begin{theorem}
\label{BEarea}
Let ${\mathbf E}:= [0,\infty) \times[0,\infty) \setminus\{(0,0)\}
$. There exist measures
$\mu$ and $\nu$ on the Borel $\sigma$-algebra on ${\mathbf E}$
such that for any Borel set $B \subseteq\mathbf E $
we have
%
\begin{eqnarray}
\label{eqnHack} \sqrt{n} \P \biggl[ \frac{ ( L(0,0), (\#C(0,0))^{2/3}
)}{n} \in B \biggr] &\to&\mu(B),
\\
\label{eqnMaxExp} \sqrt{n} \P \biggl[ \frac{ ( L(0,0), (D_{\max}(0,0))^{1/2}
 )}{n} \in B \biggr] &\to&\nu(B),
\end{eqnarray}
with $\mu$ and $\nu$ being given by
\begin{eqnarray*}
\mu(B) &=& \int\!\!\int_B \frac{3\sqrt{v}}{4 \sqrt{2\pi}\gamma
_0^2pt^3}f_{I^+_0}
\biggl(\frac{v^{3/2}}{
\gamma_0 p\sqrt{2t^3}}\biggr)\,dv \,dt,
\\
\nu(B) &=& \int\!\!\int_B \frac{v}{\sqrt{2\pi}\gamma
_0^2pt^2}f_{M^+_0}
\biggl(\frac{v^2}{\gamma_0 p\sqrt{2t}}\biggr) \,dv \,dt
\end{eqnarray*}
where $f_{I^+_0}$ and $f_{M^+_0}$ denote the density functions of
$I^+_0$ and
$M^+_0$, respectively.
Moreover, for $\lambda, \tau> 0$, we have
%
\begin{eqnarray}\label{eqnTrivialCasesHack}
&& \sqrt{n} \P \biggl[ \frac{ ( L(0,0), (\#C(0,0))^{\alpha}
)}{n} \in(\tau, \infty) \times(
\lambda, \infty) \biggr]
\nonumber\\[-8pt]\\[-8pt]\nonumber
&&\qquad = \cases{ 0, &\quad if $\alpha< \displaystyle\frac{2}{3}$,
\vspace*{4pt}\cr
\displaystyle\frac{1}{\sqrt{\pi\tau\gamma^2_0}}, &\quad if $\alpha> \displaystyle\frac{2}{3}$}
\end{eqnarray}
and
%
\begin{eqnarray}
\label{eqnTrivialCasesHurst}
&& \sqrt{n} \P \biggl[ \frac{ ( L(0,0), (D_{\max} (0,0))^{\alpha
} )}{n} \in(\tau, \infty) \times(
\lambda, \infty) \biggr]
\nonumber\\[-8pt]\\[-8pt]\nonumber
&&\qquad = \cases{ 0, &\quad if $\displaystyle \alpha< \frac{1}{2}$,
\vspace*{4pt}\cr
\displaystyle \frac{1}{\sqrt{\pi\tau\gamma^2_0}}, &\quad if $\displaystyle\alpha> \frac{1}{2}$.}
\end{eqnarray}
\end{theorem}

The estimates of the densities $ f_{I^+_0} $ and $ f_{M^+_0} $ imply
that $ \mu$ and $\nu$ are
finite measures on ${\mathbf E}$.
An immediate consequence of the above theorem is the following.
%
\begin{cor}
\label{corBEarea}
As $n\rightarrow\infty$ for $ u > 0 $, we have:
\begin{longlist}[(ii)]
\item[(i)]
$\sqrt{n}\mathbb{P} ( \# C (0,0) > \sqrt{2n^3}\gamma_0pu )
\to\frac{1}{2\sqrt{\pi}\gamma_0}\int_{0}^{ \infty} t^{- \sfrac
{3}{2}} \wb{F}_{I^{+}_{0}}
(ut^{- \sfrac{3}{2}}) \,dt $,\vspace*{2pt}
\item[(ii)] $\sqrt{n}\mathbb{P} ( D_{\max}(0,0) > \sqrt{2n}\gamma
_0pu )
\to\frac{1}{2\sqrt{\pi}\gamma_0}\int_{0}^{\infty} t^{ - \sfrac
{3}{2}} \wb{F}_{M^{+}_{0}}
(ut^{- \sfrac{1}{2}}) \,dt $,
\end{longlist}
where $F_{I^{+}_{0}}$ and $F_{M^{+}_{0}}$ are the distribution
functions of
$I^{+}_{0}$ and $M^{+}_{0}$, respectively, and $\wb{F}_{I^{+}_{0}}:=
1 - F_{I^{+}_{0}}$,
$\wb{F}_{M^{+}_{0}}:= 1 - F_{M^{+}_{0}}$.
\end{cor}
%

The proofs of the above theorems are based on a scaling of the process.
In the next section, we define a dual graph and show that as processes,
under a suitable scaling, the original and the dual processes
converge jointly to the Brownian web and its dual in distribution (the
double Brownian web). This
invariance principle is used in Sections~\ref{ClusterThm} and \ref{BMBEarea}
to prove the theorems.
In this connection, it is worth noting that in Proposition \ref
{propDualwedgealternate}, we have
provided an alternate characterization of the dual of Brownian web
which is of independent interest.
This characterization is suitable for proving the joint convergence of
coalescing noncrossing path family and its dual to the double Brownian
web and has been used in
Theorem \ref{theoremGRSDual-DobleBW} to achieve the required convergence.

We should mention here that the Brownian web appears as a universal
scaling limit
for various network models (see Fontes~et al. \cite{FINR04},
Ferrari, Fontes and Wu \cite{FFW05},
Coletti, Fontes and Dias \cite{CFD09}).
It is reasonable to expect that with suitable modifications our method
will give similar results in other network models. Our results will
hold for any network model which admits a dual and satisfies (i)
conditions listed in Remark \ref{remDualGraph},
(ii) the scaled model and its dual converges weakly to the double
Brownian web (see Section~\ref{SecDual}) and (iii) a certain sequence
of counting random variables are uniformly integrable (see Lemma \ref
{lemUI}). In this sense, our result can be considered as a
universality class result.

\section{Dual process and the double Brownian web}\label{SecDual}
\subsection{Dual process}\label{SubSecDualProcess}

\begin{figure}[b]

\includegraphics{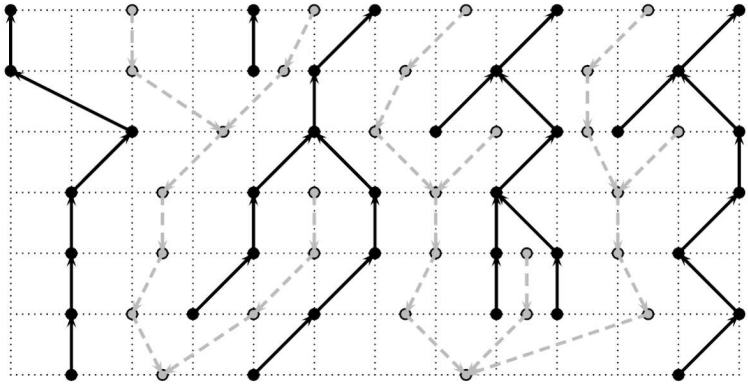}

\caption{The black points are open vertices, the gray points are the
vertices of the dual
process and the gray (dashed) paths are the dual paths.}\label{GRSDual}
\end{figure}

For the graph $\mathcal G$, we now describe a dual process such that
the set of ancestors $ C(x,t)$
(defined in the previous section) of a vertex $ (x,t) \in V$ is bounded
by two dual paths.
The dependency inherent in the graph $\mathcal G$ implies that,
although the cluster is bounded by two dual paths,
these paths are \textit{not given by independent random walks}.
The dual vertices are precisely the mid-points between two consecutive
open vertices
on each horizontal line $\{y = n\}, n \in{\mathbb Z}$ with each
dual vertex having a unique offspring dual vertex in the negative
direction of the $y$-axis.
Before giving a formal definition, we direct the attention of the
reader to Figure~\ref{GRSDual}.

%
%
For $ \bu \in{\mathbb Z}^2$, we define
%
\begin{eqnarray}
\label{DualIncrement} J^{+}_{\bu} &:=& \inf\bigl\{k: k \geq1,
\bigl( \bu(1)+k, \bu(2) \bigr) \in V\bigr\},
\nonumber\\[-8pt]\\[-8pt]\nonumber
J^{-}_{\bu} &:=& \inf\bigl\{k: k\geq1, \bigl( \bu(1) -k,
\bu(2)\bigr) \in V\bigr\}.
\end{eqnarray}
%
Next, we define $r(\bu):= ( \bu(1)+ J^{+}_{\bu},\bu(2))$ and $l(\bu
):= (\bu(1) - J^{-}_{\bu},\bu(2))$, as the first
open point to the right (\textit{open right neighbour}) and the first open
point to the left (\textit{open left neighbour})
of $\bu$, respectively.
For $(x,t)\in V$, let $\hat{r}(x,t):= (x + J^{+}_{(x,t)}/2,t)$
and $\hat{l}(x,t):= (x -J^{-}_{(x,t)}/2,t)$
denote,\vspace*{1pt} respectively, the right dual neighbour and the left dual
neighbour of $(x,t)$ in the dual vertex set.
Finally, the dual vertex set is given by
\[
\wh{V}:= \bigl\{\hat{r}(x,t),\hat{l}(x,t):(x,t) \in V \bigr\}.
\]
For a vertex $ (u,s)\in\wh{V} $, let $ (v,s-1) \in\wh{V}
$ be such that the straight line segment
joining $ (u,s)$ and $(v,s-1)$ does not cross any edge in $ {\mathcal
G}$. The dual edges are edges joining all
such $ (u,s)$ and $ (v,s-1) $. Formally, for $(u,s)\in\wh{V}$,
we define
%
\begin{eqnarray}
\label{eqndefalar} a^{l} (u,s) &:=& \sup\bigl\{ z: (z,s-1)\in V,
h(z,s-1) (1) < u \bigr\},
\nonumber\\[-8pt]\\[-8pt]\nonumber
a^{r} (u,s) &:=& \inf\bigl\{z: (z,s-1) \in V, h(z,s-1) (1) > u \bigr\}
\end{eqnarray}
and set $ \hat{h} ( u, s):= ( (a^{l} (u,s)+a^{r} (u,s))/2, s- 1)
$. Note that $ ( a^{r} (u,s), s-1) $ and
$ (a^{l} (u,s), s-1) $ are the nearest vertices in $V$ to the right and
left, respectively,
of the dual vertex $ \hat{h} ( u, s) $.
Finally, the edge set of the dual graph $\widehat{{\mathcal G}}:=
(\wh{V}, \widehat{E})$
is given by
\[
\widehat{E}:= \bigl\{ \bigl\langle(u,s), \hat{h}(u,s) \bigr\rangle: (u,s)
\in \wh{V}\bigr\}.
\]
%
%

%
\begin{remark}
\label{remDualGraph}
Note that the vertex set of the dual graph is a subset of $\frac
{1}{2}{\mathbb Z} \times{\mathbb Z}$.
Before we proceed, we list some properties of the graph $\mathcal G$
and its dual $\widehat{{\mathcal G}}$.
\begin{longlist}[(2)]
\item[(1)] ${\mathcal G}$ uniquely
specifies the dual graph $\widehat{{\mathcal G}}$ and
the dual edges do not intersect the original edges. The construction
ensures that
$\widehat{{\mathcal G}} $ does not contain any circuit.

\item[(2)] For $ (x,t) \in V$, the cluster $C(x,t)$ is enclosed
within the dual paths starting from $\hat{r}(x,t)$ and $\hat{l}(x,t)$.
The boundedness of $C(x,t)$ for every $(x,t)\in V$ implies that
these two dual paths coalesce, thus $\widehat{{\mathcal G}} $ is a
single tree.

\item[(3)] Since paths starting from any two open vertices in the
original graph
coalesce and the dual edges do not cross the original edges, there is
no bi-infinite
path in $\widehat{{\mathcal G}} $. 
\end{longlist}
\end{remark}

We now obtain a Markov process from the dual paths.
Fix $(u,s)\in\wh{V}$ and for $k\geq1$, set $\hat{h}^k(u,s)
:= \hat{h}(\hat{h}^{k-1}(u,s))$
where $ \hat{h}^0(u,s):= (u,s)$. Letting $ \widehat
{X}^{(u,s)}_{k}$ denote the first coordinate of $\hat{h}^k(u,s)$,
it may be observed that
%
$\widehat{X}^{(u,s)}_{k+1}$ is a function of
$ \widehat{X}^{(u,s)}_{k}$ and the collection of random variables 
$ \{ (B_{\bu}, U_{\bu}): \bu(2) = s -k-1\in{\mathbb Z} \} $.
%
Thus, by the
random mapping representation (see, e.g.,
Levin, Peres, and Wilmer
\cite{LPW08}) we have the following.
%
\begin{prop}
\label{DualMarkov}
For $(u,s)\in\wh{V}$, the process $\{\widehat{X}^{(u,s)}_k:
k\geq0\}$ is a
time homogeneous Markov process.
\end{prop}

Before we proceed, we make the following observations about the
transition probabilities of the
Markov process. Let $ G $ be a geometric random variable taking values
in $ \{1, 2, \dotsc\}$, that is,
$ \P( G = l ) = p (1-p)^{l-1} $ for $ l \geq1 $. For any $\bu\in
{\mathbb Z} \times{\mathbb Z}$, the random
variables $J^{+}_{\bu}$ and $J^{-}_{\bu}$ are i.i.d.\vspace*{1pt} copies of the
geometric random variable $ G$
independent of $ B_{ \bu} $.
Further, if $ \bu_1, \bu_2 \in{\mathbb Z}^2 $ are such that $ \bu
_1(1) \geq\bu_2 (1) - 1 $ and $ \bu_1 (2)
= \bu_2 (2) $, the random variables $J^{+}_{\bu_1}$
and $ J^{-}_{\bu_2}$ are also independent. Now, for $ u \notin
{\mathbb Z} $ and, $ v \in{\mathbb Z}/2 $, we have
%
\begin{eqnarray}
\label{eqntranprobnonint} \P\bigl( \widehat{X}^{(u,s)}_1 -
\widehat{X}^{(u,s)}_0 = v \mid \widehat{X}^{(u,s)}_0
= u \bigr) & =& \P\bigl( J^{+}_{(u-1/2,s-1)} - J^{-}_{(u+1/2,s-1)}
= 2v \bigr)
\nonumber\\[-8pt]\\[-8pt]\nonumber
& =& \P( G_1 - G_2 = 2v ),
\end{eqnarray}
where $ G_1 $ and $ G_2 $ are i.i.d. copies of $G$, defined above. If $
u \in{\mathbb Z} $ and $ v \in{\mathbb Z}/2 $,
we have, using notations from above
%
\begin{eqnarray}
\label{eqntranprobint}
&& \P\bigl( \widehat{X}^{(u,s)}_1 -
\widehat{X}^{(u,s)}_0 = v \mid \widehat{X}^{(u,s)}_0
= u \bigr)
\nonumber\\[-8pt]\\[-8pt]\nonumber
&&\qquad
 = (1-p) \P(
G_1 - G_2 = 2v ) + p \P(G = 2v)/2 + p \P(G = - 2v)/2,
\end{eqnarray}
where $ G_1 $ and $ G_2 $ are as above. It is therefore obvious that
the transition
probabilities of $\widehat{X}^{(u,s)}_k$ depend on
whether the present state is an integer or not.

From equations (\ref{eqntranprobnonint}) and (\ref
{eqntranprobint}), we state the following.
%
\begin{prop}
\label{propMartingale}
For any $(u,s) \in\wh{V}$, $ \{ \widehat{X}^{(u,s)}_k: k \geq0
\}$
is an $L^2$-martingale with respect to the filtration $ {\mathcal F}_k
:= \sigma
( \{ B_{\bu}, U_{ \bu}: \bu\in{\mathbb Z}^2, \bu(2) \geq s - k \})$.
\end{prop}

\subsection{Dual Brownian web}\label{SubSecDualBW}


In this section, we briefly describe
the dual Brownian web $\widehat{{\mathcal W}}$ associated with
${\mathcal W}$ and present an
alternate characterization of the dual Brownian web $\widehat
{{\mathcal W}}$.

The Brownian web
(studied extensively by Arratia \cite{A79,A81},
T{\'o}th and Werner \cite{TW98},
Fontes~et al. \cite{FINR04})
may be viewed as a collection
of one-dimensional coalescing Brownian motions starting from every
point in the space time
plane $\R^2$.
We recall relevant details from Fontes~et al. \cite{FINR04}.

Let $\R^{2}_c$ denote the completion of the space time plane $\R^2$ with
respect to the metric
\[
\rho\bigl((x_1,t_1),(x_2,t_2)
\bigr):= \bigl\llvert \tanh(t_1)-\tanh(t_2)\bigr\rrvert
\vee\biggl\llvert \frac{\tanh(x_1)}{1+\llvert   t_1\rrvert   } -\frac{\tanh(x_2)}{1+\llvert   t_2\rrvert   } \biggr\rrvert.
\]
%
As a topological space $\R^{2}_c$ can be identified with the
continuous image of $[-\infty,\infty]^2$ under a map that identifies
the line
$[-\infty,\infty]\times\{\infty\}$ with the point $(\ast,\infty
)$, and the line
$[-\infty,\infty]\times\{-\infty\}$ with the point $(\ast,-\infty)$.
A path $\pi$ in $\R^{2}_c$ with starting time $\sigma_{\pi}\in
[-\infty,\infty]$
is a mapping $\pi:[\sigma_{\pi},\infty]\rightarrow[-\infty,\infty
] \cup\{ \ast\}$ such that
$\pi(\infty)= \ast$ and, when $\sigma_\pi= -\infty$, $\pi
(-\infty)= \ast$.
Also $t \mapsto(\pi(t),t)$ is a continuous
map from $[\sigma_{\pi},\infty]$ to $(\R^{2}_c,\rho)$.
We then define $\Pi$ to be the space of all paths in $\R^{2}_c$ with
all possible starting times in $[-\infty,\infty]$.
The following metric, for $\pi_1,\pi_2\in\Pi$
\begin{eqnarray*}
d_{\Pi} (\pi_1,\pi_2)&:=& \max \biggl\{\bigl
\llvert \tanh(\sigma_{\pi
_1})-\tanh(\sigma_{\pi_2})\bigr\rrvert,
\\
&&{}\sup_{t\geq\sigma_{\pi_1}\wedge
\sigma_{\pi_2}} \biggl\llvert \frac{\tanh(\pi_1(t\vee\sigma_{\pi
_1}))}{1+\llvert   t\rrvert   }-
\frac{
\tanh(\pi_2(t\vee\sigma_{\pi_2}))}{1+\llvert   t\rrvert   }\biggr\rrvert \biggr\}
\end{eqnarray*}
%
makes $\Pi$ a complete, separable metric space.

%
\begin{remark}
\label{remMetricConv}
Convergence in this metric can be described as locally uniform
convergence of paths as
well as convergence of starting times. Therefore, for any $ \varepsilon>
0$ and $m>0$,
we can choose $ \varepsilon_1 ( = f ( \varepsilon, m)) > 0 $ such that for
$ \pi_1, \pi_2
\in\Pi$ with $\{(\pi_i(t),t):t\in[\sigma_{\pi_i},m]\}\subseteq
[-m,m]\times[-m,m]$ for $i=1,2$,
$ d_{\Pi} ( \pi_1, \pi_2 ) < \varepsilon_1 $ implies that $ \Vert  ( \pi
_1 (\sigma_{\pi_1}),
\sigma_{\pi_1} ) - ( \pi_2 (\sigma_{\pi_2}), \sigma_{\pi_2} )
\Vert _2 < \varepsilon$ and $ \sup\{
\llvert    \pi_1 (t) - \pi_2 (t) \rrvert: t \in[\max\{\sigma_{\pi_1},\sigma
_{\pi_2}\},m] \} < \varepsilon$.
We will use this later several times. 
\end{remark}

Let ${\mathcal H}$ be the space of compact subsets of $(\Pi,d_{\Pi})$
equipped with
the Hausdorff metric $d_{{\mathcal H}}$.
The Brownian web ${\mathcal W}$ is a random
variable taking values in the complete separable metric space
$({\mathcal H},d_{{\mathcal H}})$.

Before introducing the dual Brownian web, we require a similar metric space
on the collection of backward paths.
As in the definition of $\Pi$, let $ \wh{\Pi}$ be the
collection of all paths $ \hat{\pi}$
with starting time $\sigma_{\hat{\pi}} \in[-\infty,\infty]$
such that
$\hat{\pi}: [-\infty, \sigma_{\hat{\pi}}] \to [-\infty
,\infty] \cup\{\ast\}$ with
$\hat{\pi} (-\infty)= \ast$ and, when $\sigma_{\hat{\pi
}} = +\infty$, $\hat{\pi}(\infty)= \ast$.
As earlier $t \mapsto(\hat{\pi}(t),t)$ is a continuous
map from $[-\infty, \sigma_{\hat{\pi}} ]$ to $(\R^{2}_c,\rho)$.
We equip $ \wh{\Pi}$ with the metric
\begin{eqnarray*}
d_{\wh{\Pi}} (\hat{\pi}_1,\hat{\pi}_2)&:=& \max \biggl\{\bigl\llvert \tanh(\sigma_{\hat{\pi}_1})-\tanh(
\sigma_{\hat{\pi}_2})\bigr\rrvert,
\\
&&{}\sup_{t\leq\sigma_{\hat{\pi}_1}\vee
\sigma_{\hat{\pi}_2}} \biggl\llvert \frac{\tanh(\hat{\pi
}_1(t\wedge\sigma_{\hat{\pi}_1}))}{1+\llvert   t\rrvert   }-
\frac{\tanh
(\hat{\pi}_2(t\wedge\sigma_{\hat{\pi}_2}))}{1+\llvert   t\rrvert   }\biggr\rrvert \biggr\}
\end{eqnarray*}
%
making $(\wh{\Pi}, d_{\wh{\Pi}})$ a complete, separable
metric space.
The complete separable metric space of compact sets of paths of
$\wh{\Pi}$ is denoted
by $(\widehat{{\mathcal H}}, d_{\widehat{{\mathcal H}}})$, where
$d_{\widehat{{\mathcal H}}}$
is the Hausdorff metric on $\widehat{{\mathcal H}}$, and let $
{\mathcal B}_{\widehat{{\mathcal H}}}$
be the corresponding Borel $\sigma$ field.

\subsection{Properties of $(\mathcal{W},\widehat{\mathcal{W}})$}

The Brownian web and its dual $( {\mathcal W},\widehat{{\mathcal W}})$
is a $({\mathcal H}\times
\widehat{{\mathcal H}}, {\mathcal B}_{{\mathcal H}}\times{\mathcal
B}_{\widehat{{\mathcal H}}})$ valued
random variable such that ${\mathcal W}$ and $\widehat{{\mathcal W}}$
uniquely determine
each other almost surely with $\widehat{{\mathcal W}}$ being equally
distributed as $-{\mathcal W}$,
the Brownian web rotated 180\tsup{o} about the origin. The interaction
between the paths
in ${\mathcal W}$ and $\widehat{\mathcal W}$ is that of Skorohod
reflection (see Soucaliuc, T{\'o}th and Werner
\cite{STW00}).

We introduce some notation to study the sets $ \{ \pi(t+s): \pi\in
{\mathcal W}, \sigma_{\pi} \leq t \} $ and
$ \{ \hat{\pi} (t - s): \hat{\pi} \in\widehat{\mathcal
W}, \sigma_{\hat{\pi}} \geq t \} $.
For a $({\mathcal H},B_{{\mathcal H}})$ valued random variable $K$ and
$t\in\R$,
let $K^{t-}:= \{\pi:\pi\in K$ and $\sigma_{\pi}\leq t\}$.
Similarly, for
a $(\widehat{{\mathcal H}},B_{\widehat{{\mathcal H}}})$ valued random
variable $\widehat{K}$
and $t\in\R$, let $\widehat{K}^{t+}:= \{\hat{\pi}: \hat{\pi}
\in\widehat{K}$ and $\sigma_{\hat{\pi}} \geq t\}$.
For $t_1, t_2 \in\R$, $t_2 > t_1$ and a $({\mathcal H},B_{{\mathcal
H}})$ valued random variable $K$,
define
%
\begin{eqnarray}
\label{eqnDefNK} {\mathcal M}_{K}(t_1,t_2) &:=& \bigl\{\pi(t_2):\pi\in K^{t_1 -}, \pi (t_2)
\in[0,1]\bigr\};
\nonumber\\[-8pt]\\[-8pt]\nonumber
\xi_{K}(t_1,t_2) &:=&\#{\mathcal
M}_{K}(t_1,t_2),
\end{eqnarray}
that is, $\xi_{K}(t_1,t_2)$ denotes
the number of distinct points in $[0,1]\times t_2$ which are on some
path in $K^{t_1-}$.
We note that for $t > 0$, ${\mathcal M}_{\mathcal W}(t_0,t_0 + t) =
{\mathcal N}_{\mathcal W}(t_0,t;0,1)$
as defined in Sun and Swart \cite{SS08}.
It is known that for all $t > 0$
the random variable $\xi_{{\mathcal W}}(t_0,t_0 + t)$ is finite almost surely
(see $(E_1)$ in Theorem 1.3 in Sun and Swart \cite{SS08}) with
%
\begin{equation}
\label{expEta} \E\bigl(\xi_{{\mathcal W}}(t_0,t_0 + t)
\bigr) = \frac{1}{\sqrt{\pi t}}.
\end{equation}
Moreover, from the known properties of $ ({\mathcal W},\widehat
{{\mathcal W}})$ the proof of the following proposition is
straightforward (for details, see Roy, Saha and Sarkar \cite{RSS15}).
%
\begin{prop}
\label{lemPropEtaPtset}
For any $t_0 < t_1$, almost surely we have:
\begin{longlist}[(iii)]
\item[(i)] ${\mathcal M}_{{\mathcal W}}(t_0,t_1) \cap\Q= \varnothing$;
%
\item[(ii)] each point in ${\mathcal M}_{{\mathcal W}}(t_0,t_1)$ is of
type $(1,1)$;
\item[(iii)] for each $x \in{\mathcal M}_{{\mathcal W}}(t_0,t_1)$,
there exists $\pi_1,
\pi_2 \in{\mathcal W}$ with $\sigma_{\pi_1} < t_0$, $\sigma_{\pi
_2} > t_0$ and $\pi_1(t_1)
= \pi_2(t_1)= x$;
%
\item[(iv)] for each $x \in{\mathcal M}_{{\mathcal W}}(t_0,t_1)$,
there exist exactly
two paths $\hat{\pi}^{(x,t_1)}_r$ and $\hat{\pi
}^{(x,t_1)}_l$ in $\widehat{{\mathcal W}}$
starting from $(x,t_1)$ with $\hat{\pi}^{(x,t_1)}_r(t) > \hat{\pi}^{(x,t_1)}_l(t)$ for all $[t_0,t_1)$.
\end{longlist}
\end{prop}

There are several ways to construct
$\widehat{{\mathcal W}}$ from ${\mathcal W}$. In this paper, we follow
the \textit{wedge characterization}
provided by Sun and Swart \cite{SS08}.
For $\pi^r,\pi^l \in{\mathcal W}$ with coalescing time $t^{\pi
^r,\pi^l}$
and $\pi^r(\max\{\sigma_{\pi^r},\sigma_{\pi^l}\})>
\pi^l(\max\{\sigma_{\pi^r},\sigma_{\pi^l}\})$,
the wedge with right boundary $\pi^r$ and left boundary $\pi^l$,
is an open set in $\R^2$ given by
%
\begin{eqnarray}\label{defwedgenoncpt}
A &=& A\bigl(\pi^r, \pi^l\bigr)
\nonumber\\[-8pt]\\[-8pt]\nonumber
&:=& \bigl
\{(y,s):\max\{\sigma_{\pi^l},\sigma _{\pi^r}\} < s <
t^{\pi^r,\pi^l}, \pi^l(s) < y <\pi^r(s) \bigr\}.
\end{eqnarray}
A path $\hat{\pi} \in\wh{\Pi}$, is said to \textit{enter
the wedge
$A$ from outside} if there exist $t_1$ and $t_2$ with $ \sigma
_{\hat{\pi}} > t_1 > t_2$ such that
$(\hat{\pi}(t_1), t_1) \notin\wb{A}$ and $ (\hat{\pi
}(t_2), t_2) \in A$.

From Theorem 1.9 in Sun and Swart \cite{SS08}, it follows that the dual Brownian
web $\widehat{\mathcal W}$ associated with
the Brownian web ${\mathcal W}$ satisfies the following wedge characterization.
%
\begin{theorem}
\label{theoremDualwedge}
Let $({\mathcal W}, \widehat{{\mathcal W}})$ be a Brownian web and its
dual. Then almost surely
\[
\widehat{{\mathcal W}} = \{\hat{\pi}: \hat{\pi} \in \wh{\Pi}
\mbox{ and does not enter any wedge in } {\mathcal W} \mbox{ from outside} \}.
\]
\end{theorem}

Because of Theorem \ref{theoremDualwedge}, for a $({\mathcal H}\times
\widehat{\mathcal H}, {\mathcal B}_{{\mathcal H}}
\times{\mathcal B}_{\widehat{\mathcal H}})$ valued random variable
$({\mathcal W}, {\mathcal Z})$ to show
that ${\mathcal Z}= \widehat{\mathcal W}$, it suffices to check that
${\mathcal Z}$ satisfies the wedge condition. Here we present an
alternate condition which
is easier to check.

%
\begin{prop}
\label{propDualwedgealternate}
Let $({\mathcal W}, {\mathcal Z})$ be a $({\mathcal H}\times\widehat
{\mathcal H}, {\mathcal B}_{{\mathcal H}}
\times{\mathcal B}_{\widehat{\mathcal H}})$ valued random variable
such that:
\begin{longlist}[(2)]
\item[(1)] for any deterministic $ (x,t) \in\R^2$, there exists a
path $\hat{\pi}^{(x,t)}
\in{\mathcal Z}$ starting at $(x,t)$ and going backward in time almost surely;

\item[(2)] paths in ${\mathcal Z}$ do not cross paths in ${\mathcal
W}$ almost surely, that is, there
does not exist any $\pi\in{\mathcal W}$, $\hat{\pi} \in
{\mathcal Z}$ and $t_1,t_2 \in
(\sigma_{\pi}, \sigma_{\hat{\pi}})$ such that $(\hat{\pi
}(t_1)-\pi(t_1))
(\hat{\pi}(t_2)-\pi(t_2))< 0$ almost surely;

\item[(3)] paths in ${\mathcal Z}$ and paths in ${\mathcal W}$ do
not coincide over any time interval
almost surely, that is, for any $ \pi\in{\mathcal W}$
and $ \hat{\pi} \in{\mathcal Z}$ and for no pair of points $ t_1
< t_2 $ with $ \sigma_{\pi}
\leq t_1 < t_2 \leq\sigma_{\hat{\pi}} $
we have $ \hat{\pi}(t) = {\pi}(t) $ for all $ t \in[t_1,t_2]$
almost surely.
\end{longlist}
Then ${\mathcal Z} = \widehat{{\mathcal W}}$ almost surely.
\end{prop}

\begin{pf}
From conditions (2) and (3), we have that $\hat{\pi} \in
{\mathcal Z}$ does not enter any wedge in ${\mathcal W}$ from outside.
Hence, ${\mathcal Z} \subseteq\widehat{\mathcal W}$.
The argument for $\widehat{\mathcal W} \subseteq{\mathcal Z}$ follows from
the \textit{fish-trap} technique introduced in the proof of\vspace*{2pt} Lemma 4.7
of Sun and Swart \cite{SS08}.
It shows that
$\widehat{\mathcal W} \subseteq\widetilde{\mathcal Z}$ almost surely for any
$({\mathcal H}, {\mathcal B}_{{\mathcal H}})$ valued random variable
$\widetilde{\mathcal Z}$ satisfying
(i) paths in $\widetilde{\mathcal Z}$ do not\vspace*{1pt} cross paths is ${\mathcal W}$
and (ii) for any deterministic
countable dense set, there exist paths in $\widetilde{\mathcal Z}$
starting from every point of that dense set (for details, see Roy, Saha and Sarkar \cite{RSS15}).
\end{pf}

\subsection{Convergence to the double Brownian web}

For any\vspace*{1pt} $(x,t) \in V$, the path $\pi^{(x,t)}$ in the random graph
$\mathcal G$ is obtained as the
piecewise linear function $\pi^{(x,t)}: [t, \infty) \to\mathbb R$
with $\pi^{(x,t)}(t+k) = h^{k}(x,t)(1)$ for every $k \geq0$ and
$\pi^{(x,t)}$ being linear in the interval $[t+k,t+k + 1]$.
Similarly, for $(x,t)\in\wh{V}$,
the dual path $\hat{\pi}^{(x,t)}$ is the piecewise linear
function $\hat{\pi}^{(x,t)}: (-\infty, t] \to\mathbb R$
with $\hat{\pi}^{(x,t)}(t-k) = \hat{h}^{k}(x,t)(1)$ for
every $k \geq0$ and
$\hat{\pi}^{(x,t)}$ being linear in the interval $[t-k-1,t-k]$.
Let ${\mathcal X}:= \{\pi^{(x,t)}:(x,t) \in V\}$ and
$\widehat{{\mathcal X}}:= \{\hat{\pi}^{(x,t)}: (x,t)\in
\wh{V}\}$ be the collection of all possible paths and dual paths
admitted by $\mathcal G$ and $\widehat{\mathcal G}$.

For a given $ \gamma> 0$ and a path $\pi$ with starting time $
\sigma_{\pi}$, the scaled path
$ \pi_n(\gamma): [ \sigma_{\pi}/n, \infty] \to
[-\infty, \infty]$ is given by $\pi_n(\gamma) (t)=
\pi(n t)/ (\sqrt{n} \gamma)$ for each $n \geq1$. Thus, the starting
time of the scaled path $\pi_n(\gamma)$ is
$\sigma_{\pi_n(\gamma) }= \sigma_{\pi}/ n $.
Similarly, for the backward path $ \hat{\pi} $, the scaled
version is
$ \hat{\pi}_n(\gamma): [ -\infty, \sigma_{\hat{\pi}}/n]
\to
[-\infty, \infty]$ given by $\hat{\pi}_n(\gamma) (t)=
\hat{\pi}(n t)/ (\sqrt{n} \gamma)$ for each $n \geq1$.
For each $ n \geq1 $, let $ {\mathcal X}_n = {\mathcal X}_n(\gamma):=
\{\pi_n^{(x,t)}(\gamma):(x,t) \in V\}$ and $\widehat{{\mathcal X}}_n
= \widehat{{\mathcal X}}_n(\gamma):=
\{ \hat{\pi}_n^{(x,t)}(\gamma):(x,t) \in\wh{V} \}$ be the
collections of all the $n$th order diffusively scaled paths and dual paths, respectively.

The\vspace*{1pt} closure $\wb{\mathcal X}_n(\gamma)$ of
${\mathcal X}_n(\gamma)$ in $(\Pi,d_{\Pi})$ and the closure
$\wb{\widehat{{\mathcal X}}}_n(\gamma)$
of $\wh{\mathcal X}_n(\gamma)$ in $(\wh{\Pi}, d_{\widehat
{\Pi}})$ are
$({\mathcal H},{\mathcal B}_{{\mathcal H}})$ and $ (\widehat{{\mathcal
H}}, {\mathcal B}_{\widehat{{\mathcal H}}})$
valued random variables, respectively. Coletti, Fontes and Dias \cite{CFD09} showed the
following.%
\begin{theorem}
\label{theoremGRS-BW}
For $\gamma_0:= \gamma_0 (p)$ as in Theorem \ref{clusterheight}, as
$n\rightarrow\infty$,
$\wb{{\mathcal X}}_n(\gamma_0)$ converges weakly to the
standard Brownian web ${\mathcal W}$.
\end{theorem}

Our main result is the joint invariance principle
for $\{(\wb{\mathcal X}_n(\gamma_0),\wb{\widehat
{{\mathcal X}}}_n(\gamma_0)): n \geq1\}$ considered as $({\mathcal
H} \times\widehat{{\mathcal H}}, {\mathcal B}_{{\mathcal H}} \times
{\mathcal B}_{\widehat{{\mathcal H}}})$ valued random variables.
%
\begin{theorem}
\label{theoremGRSDual-DobleBW}
$\{ (\wb{{\mathcal X}}_n(\gamma_0), \wb{\widehat{{\mathcal X}}}_n(\gamma_0)): n \geq0\}$ converges weakly
to $({\mathcal W}, \widehat{{\mathcal W}})$ as $n\rightarrow\infty$.
\end{theorem}

We require the following propositions to prove Theorem \ref
{theoremGRSDual-DobleBW}. We say that
$ \{ \widehat{W}^{(x,t)}(u): u \leq t \}$ is a Brownian motion going
\textit{back in time} if
$ \widehat{W}^{(x,t)}(t-s): = W(t+s),   s \geq0$ where $\{W(u): u
\geq t\}$ is a
Brownian motion with $W(t) = x$.
%
\begin{prop}
\label{propDual-Bmotion}
For any deterministic point $(x,t)\in\R^2$, there exists a sequence
of paths $\hat{\theta}^{(x,t)}_n \in\widehat{{\mathcal
X}}_n(\gamma_0)$
which converges in distribution to $ \widehat{W}^{(x,t)}$.
\end{prop}

\begin{pf}
For any $(x,t) \in\R^2$ fix $ t_n
= \lfloor n t \rfloor$ and $ x_n = \max\{ \lfloor\sqrt{n} \gamma_0
x \rfloor+ j: j \leq0,
( \lfloor\sqrt{n} \gamma_0 x \rfloor+ j, t_n ) \in\wh{V} \}
$. Let $ \hat{\theta}^{(x,t)}_n \in
\widehat{{\mathcal X}}_n (\gamma_0) $ be the scaling of the path
\mbox{$\hat{\pi}^{(x_n,t_n)} \in\widehat{{\mathcal X}} $}.

Since $ {\mathcal G}$ is invariant under translation by lattice points
and $\widehat{\mathcal G}$ is uniquely
determined by ${\mathcal G}$, the conditional distribution of $ \{
(x_n,t_n) + \hat{h}^j(0,0): j \geq0 \} $ given
$ (0,0) \in\wh{V} $ is the same as that of $ \{ \hat{h}^j(x_n,t_n ): j \geq0 \} $.
We observe that
$ ( x_n /( \sqrt{n} \gamma_0), t_n / n ) \to(x,t)$ as $n \to\infty
$ almost surely. Hence, it suffices to prove that the scaled
dual path starting from $ (0,0)$ given $ (0,0)\in\wh{V}$
converges in distribution to $ \widehat{W}^{(0,0)}$.

From Proposition \ref{propMartingale}, we see that $ \widehat
{X}_j^{(0,0)} = \hat{h}^j(0,0) (1) $ is
an $L^2$ martingale with respect to the filtration $ \sigma( \{
B_{(z,s)}, U_{ (z,s)}: z \in{\mathbb Z}, s \geq - k \})$.
Let
\[
\eta_n ( u ):= s_n^{-1} \bigl[
\widehat{X}_j^{(0,0)} + \bigl(\widehat {X}_{j+1}^{(0,0)}
- \widehat{X}_j^{(0,0)} \bigr) \bigl(u s_{n}^2
- s_j^2\bigr) / \bigl( s_{j+1}^2 -
s_j^2\bigr) \bigr]
\]
for\vspace*{2pt} $ u \in[0,\infty) $ and $ s_j^2 \leq u s_n^2 < s_{j+1}^2 $,
where $ s_n^2 = \sum_{j=1}^n \E( (\widehat{X}_{j}^{(0,0)} - \widehat
{X}_{j-1}^{(0,0)} )^2 ) $.
We know $ \eta_n $ converges in distribution
to a standard Brownian motion (see Theorem~3,~\cite{B71}). Since $ s_n^2 / (n \gamma_0^2) \to1 $,
it can be seen that $ \sup_{u \in[0,M]} \llvert    \eta_n (u) - \hat{\theta}_n^{(0,0)} (-u) \rrvert    \to0 $ in probability
for any $ M > 0 $.
So by Slutsky's theorem, we conclude that
$\hat{\theta}_n^{(0,0)}$ converges in distribution
to a standard Brownian motion going backward in time.
\end{pf}

The next result helps in estimating the probability that a direct path
and a dual path
stay close to each other for some time period.
Given $m \in\N$ and $\varepsilon, \delta>0$, we define the event
\begin{eqnarray*}
B^{\varepsilon}_n &= & B^{\varepsilon}_n(\delta, m)
\\
&:=& \bigl\{ \mbox {there exist }\pi_1^n,
\pi_2^n, \pi_3^n \in{\mathcal
X}_n \mbox{ such that } \sigma_{\pi_1^n},
\sigma_{\pi_2^n} \leq 0,
\\
&&{} \sigma_{\pi_3^n} \leq\lfloor n\delta\rfloor/n, \pi_1^n(0)
\in [-m,m], \bigl\llvert \pi_1^n(0) -
\pi_2^n(0)\bigr\rrvert < \varepsilon,\mbox{ with}
\\
&&{} \pi_1^n\bigl(\lfloor n\delta\rfloor/n\bigr) \neq
\pi_2^n\bigl(\lfloor n\delta \rfloor/n\bigr)
\mbox{ and } \bigl\llvert \pi_1^n\bigl( \lfloor n
\delta\rfloor/n\bigr) - \pi _3^n\bigl( \lfloor n\delta
\rfloor/n\bigr)\bigr\rrvert < \varepsilon,\mbox{ with}
\\
&&{}\pi_1^n\bigl( 2 \lfloor n\delta\rfloor/n\bigr) \neq
\pi_3^n\bigl( 2 \lfloor n\delta\rfloor/n\bigr) \bigr\}.
\end{eqnarray*}

%
\begin{lemma}
\label{lemmacoalescenceGRS}
For any $m \in\N$ and $\varepsilon, \delta> 0$, we have
\[
\P\bigl( B^{\varepsilon}_n (\delta, m) \bigr) \leq C_1
(\delta, m) \varepsilon,
\]
where $ C_1 (\delta, m) $ is a positive constant, depending only on $
\delta$ and $ m $.
\end{lemma}

\begin{pf}
Let $ D^{\varepsilon}_n $ be the unscaled version of the event $
B^{\varepsilon}_n $, that is,
\begin{eqnarray*}
D^{\varepsilon}_n &:=& \bigl\{ \mbox{there exist } (x,0), (y,0),
\bigl(z, \lfloor n\delta\rfloor\bigr) \in V \mbox{ such that}
\\
&&{} x \in[ - m \sqrt{n} \gamma_0, m \sqrt{n} \gamma_0 ],
\llvert x - y \rrvert < \sqrt{n} \varepsilon\gamma_0\mbox{ and }
h^{\lfloor n\delta \rfloor} (x, 0) \neq h^{\lfloor n\delta\rfloor} (y, 0),
\\
&&{} \bigl\llvert h^{\lfloor n\delta\rfloor} (x, 0) (1) - z \bigr\rrvert < \sqrt{n} \varepsilon
\gamma_0, h^{2 \lfloor n\delta\rfloor} (x, 0) \neq h^{\lfloor n\delta\rfloor} \bigl(z,
\lfloor n\delta\rfloor \bigr) \bigr\}.
\end{eqnarray*}
%

%
\begin{figure}

\includegraphics{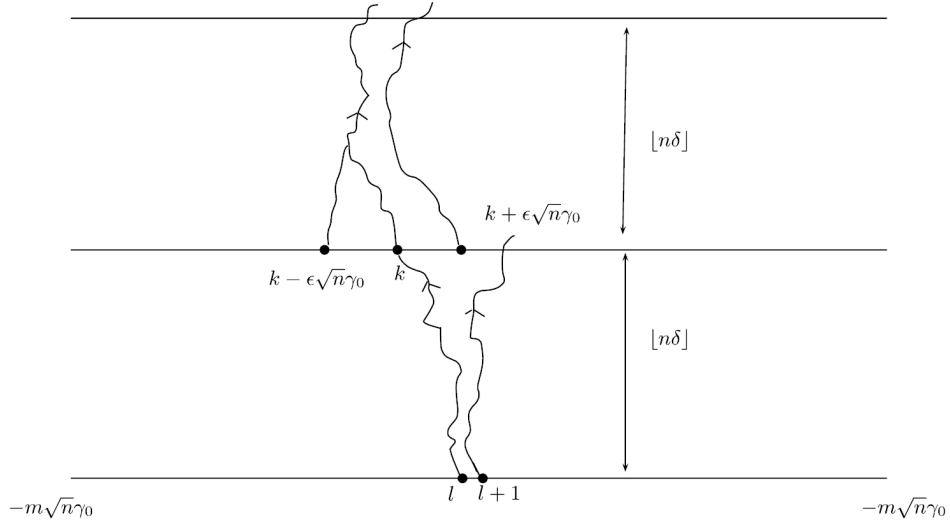}

\caption{The vertices $(l,0)$ and $(l+1,0)$ and the corresponding
vertex $(k, \lfloor n\delta\rfloor)$ as required in the proof of
Lemma \protect\ref{lemmacoalescenceGRS}.}\label{lem211}
\end{figure}

On the event $D^{\varepsilon}_n$ there exists $l \in[- m \sqrt{n}
\gamma_0, m \sqrt{n} \gamma_0]
\cap{\mathbb Z}$ such that the unscaled paths starting from $(l,0)$
and $(l+1, 0)$
(as in Figure~\ref{lem211}) do not meet in time $\lfloor n\delta
\rfloor$---an event which occurs with probability at most $C_2/\sqrt{n \delta}$
for some constant $C_2 > 0$ (see Theorem 4 of Coletti, Fontes and Dias \cite{CFD09}).
Supposing $h^{\lfloor n\delta\rfloor} (l, 0)(1) = k$, two unscaled
paths, one starting from a vertex $\lfloor\sqrt{n} \varepsilon\gamma_0
\rfloor$ distance to the left of~$k$ and the other starting from a
vertex $\lfloor\sqrt{n} \varepsilon\gamma_0 \rfloor$ distance to the
right of $k$, do not meet in time $\lfloor n\delta\rfloor$ has a
probability at most $C_2 2\sqrt{n} \varepsilon\gamma_0 /\sqrt{n \delta
}$ for all $k \in{\mathbb Z}$.
Thus, summing over all possibilities of $l$ and $k$ and using Markov
property we have
\begin{eqnarray*}
\P\bigl(D^{\varepsilon}_n\bigr) &\leq&\P\Biggl(\bigcup
_{l= - 2m \sqrt{n} \gamma_0}^{2m
\sqrt{n} \gamma_0} \bigcup_{k \in{\mathbb Z}}
\bigl\{ h^{\lfloor n\delta
\rfloor} (l, 0) (1) = k \neq h^{\lfloor n\delta\rfloor} (l+1, 0) (1) \mbox{ and}
\\
&&{} h^{\lfloor n\delta\rfloor} \bigl(k - \lfloor\sqrt{n} \varepsilon
\gamma_0 \rfloor, \lfloor n\delta\rfloor\bigr) \neq h^{\lfloor n\delta
\rfloor}
\bigl(k + \lfloor\sqrt{n} \varepsilon\gamma_0 \rfloor, \lfloor n\delta
\rfloor\bigr)\bigr\}\Biggr)
\\
&\leq&\sum_{l= - 2m \sqrt{n} \gamma_0}^{ 2m \sqrt{n} \gamma_0} \frac{2 C_2 \sqrt{n} \varepsilon\gamma_0 }{ \sqrt{n \delta}}
\sum_{k \in{\mathbb Z}} \P\bigl\{ h^{\lfloor n\delta\rfloor} (l, 0) (1) = k
\neq h^{\lfloor
n\delta\rfloor} (l+1, 0) (1)\bigr\}
\\
&\leq&\sum_{l= - 2m \sqrt{n} \gamma_0}^{ 2m \sqrt{n} \gamma_0} \frac{2 C_2 \sqrt{n} \varepsilon\gamma_0 }{ \sqrt{n \delta}} \P
\bigl\{ h^{\lfloor n\delta\rfloor} (l, 0) (1) \neq h^{\lfloor n\delta\rfloor
} (l+1, 0) (1)\bigr\}
\\
&\leq& \sum_{l= - 2m \sqrt{n} \gamma_0}^{ 2m \sqrt{n} \gamma_0} \frac{2 C_2 \sqrt{n} \varepsilon\gamma_0 }{ \sqrt{n \delta}}
\frac
{C_2}{\sqrt{n \delta}}
\\
&\leq& C_1 (\delta, m) \varepsilon.
\end{eqnarray*}\upqed
\end{pf}

\begin{pf*}{Proof of Theorem \ref{theoremGRSDual-DobleBW}}
Since $\widehat{{\mathcal X}}$ consists of noncrossing paths only, Proposition
\ref{propDual-Bmotion} implies the tightness of the family $\{\wb
{\widehat{{\mathcal X}}}_n: n\geq1\}$
(see Proposition B.2 in the Appendix of Fontes~et al. \cite{FINR04}).
The joint family $\{(\wb{{\mathcal X}}_n,\wb{\widehat{{\mathcal X}}}_n): n\geq1\}$
is tight since
each of the two marginal families 
is tight. To prove Theorem \ref{theoremGRSDual-DobleBW}, it suffices
to show that for any
subsequential limit $({\mathcal W}, {\mathcal Z})$ of $\{(\wb
{{\mathcal X}}_n,\wb{\widehat{{\mathcal X}}}_n): n\geq1\}$,
the random variable ${\mathcal Z}$ satisfies the conditions given in
Proposition \ref{propDualwedgealternate}.

Consider a convergent subsequence of $\{(\wb{\mathcal X}_n,
\wb{\widehat{{\mathcal X}}}_n): n\geq1\}$ such that $({\mathcal
W},{\mathcal Z})$ is its weak limit
and by Skorohod's representation theorem, we may
assume that the convergence happens almost surely.
For ease of notation, we denote the convergent subsequence by itself.

From Proposition \ref{propDual-Bmotion}, it follows that for any deterministic
$(x,t) \in\R^2$ there exists a path $\hat{\pi} \in{\mathcal
Z}$ starting at $(x,t)$
going backward in time almost surely.

Since $(\wb{\mathcal X}_n,
\wb{\widehat{{\mathcal X}}}_n)$ converges to $({\mathcal
W},{\mathcal Z})$ almost surely, if a dual path
in ${\mathcal Z}$ crosses a path in $\mathcal W$, there exists a dual
path in
$\wh{\mathcal X}_n$ which crosses a path in ${\mathcal X}_n$, for
some $n \geq1$,
yielding a contradiction.
Hence, the paths in ${\mathcal Z}$ do not cross paths in $\mathcal W$
almost surely
(for details, see Roy, Saha and Sarkar \cite{RSS15}).

Now, to prove that condition (3) in Proposition \ref
{propDualwedgealternate} is satisfied,
we define the following event: for $ \delta> 0 $ and positive integer
$ m \geq1 $, let
\begin{eqnarray*}
A (\delta, m ) &:=& \bigl\{\mbox{there exist paths } \pi\in {\mathcal W}
\mbox{ and }\hat{\pi} \in{\mathcal Z} \mbox{ with }
\sigma_{\pi},\sigma_{\hat{\pi}} \in (-m,m),
\\
&&{}\mbox{and there exists } t_0 \mbox{ such that }
\sigma_{\pi
} < t_0 < t_0 + \delta<
\sigma_{\hat{\pi}},
\\
&&{}\mbox{and } {-}m < \pi(t) = \hat{\pi}(t) < m \mbox { for all } t
\in[t_0, t_0+\delta] \bigr\}.
\end{eqnarray*}
It is enough to show that for any fixed $ \delta> 0 $ and for $m \geq1$,
we have $ \P ( A (\delta,  m )  ) = 0 $.

We present here the idea of the proof; more details are available in
Roy, Saha and Sarkar \cite{RSS15}.
Fix $\varepsilon> 0$. Since we are in a setup where the scaled paths
converge almost surely,
for all large $n$ there exist $\pi^n_1 \in{\mathcal X}_n$ and
$\hat{\pi}^n
\in\wh{\mathcal X}_n$ within $\varepsilon$ distance of $\pi$ and
$\hat{\pi}$, respectively.
Using the fact that a dual vertex lies in the middle of two open
vertices and the forward paths
cannot cross the dual paths, it follows that for all large $n$
there exist $\pi^n_2, \pi^n_3 \in\wh{\mathcal X}_n$
such that:
\begin{longlist}[(a)]
\item[(a)] $\max\{\llvert   \pi^n_1(\sigma_{\pi^n_2})-\pi^n_2(\sigma_{\pi^n_2})\rrvert,
\llvert   \pi^n_1(\sigma_{\pi^n_3})-\pi^n_3(\sigma_{\pi^n_3})\rrvert   \} <
4\varepsilon$;\vspace*{2pt}
\item[(b)] $\pi^n_1(\sigma_{\pi^n_2}+\delta/3) \neq\pi
^n_2(\sigma_{\pi^n_2}+ \delta/3)$ and
$\pi^n_1(\sigma_{\pi^n_3}+\delta/3) \neq\pi^n_3(\sigma_{\pi
^n_3}+ \delta/3)$.
\end{longlist}
This gives us that
$A (\delta, m ) \subseteq\liminf_{n \to\infty} \bigcup_{j =
1}^{\lfloor6m/\delta\rfloor} B^{4\varepsilon}_n( \delta/3, 2m; j ) $.

%
\begin{figure}

\includegraphics{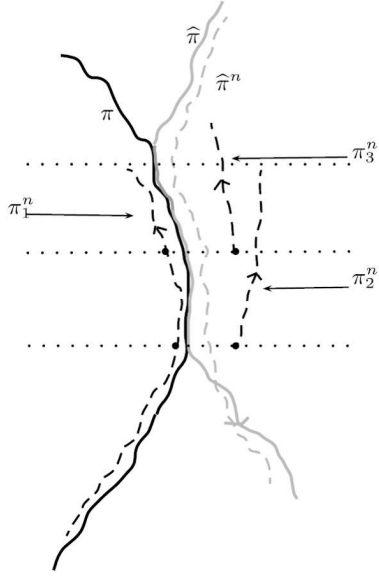}

\caption{The event $A(\delta,m)$. The bold paths are from $({\mathcal
W}, \widehat{{\mathcal W}})$
and the approximating dashed paths are from $({\mathcal X}_n, \widehat
{\mathcal X}_n)$.}\label{figDoubleBweb}
\end{figure}

Here, $B^{4\varepsilon}_n( \delta/3, 2m; j )$ is a translation
of the event $ B^{4\varepsilon}_n(\delta/3, 2m) $, considered in Lemma
\ref{lemmacoalescenceGRS};
translated such that the starting time of the paths $\pi^1_n$ and $\pi
^2_n$ are shifted by $ -m +
j \lfloor n \delta/ 3 \rfloor/ n $ (see Figure~\ref{figDoubleBweb}).

By\vspace*{1pt} translation invariance of our model and Lemma \ref
{lemmacoalescenceGRS}, for all $ n \geq1$
we have $ \P( B^{4\varepsilon}_n(\delta/3, 2m; j )) \leq4 C_1 ( \delta
/3, 2m) \varepsilon$.
This completes the proof.
\end{pf*}

\section{Proof of Theorem \texorpdfstring{\protect\ref{clusterheight}}{1.2}}\label{ClusterThm}

%
Let $\xi:= \xi_{{\mathcal W}}(0,1)$ and $\xi_n:=
\xi_{\wb{\mathcal X}_n}(0,1)$ be as defined in (\ref{eqnDefNK}).
The proof of Theorem \ref{clusterheight} follows from the following
proposition.

%
\begin{prop}
\label{propWeakConv1}
$\E[\xi_n] \to\E[\xi]$ as $n \to\infty$.
\end{prop}

We first complete the proof of Theorem \ref{clusterheight} assuming
Proposition \ref{propWeakConv1}.

\begin{pf*}{Proof of Theorem \ref{clusterheight}}
Using the translation
invariance of our model, we have
\begin{eqnarray*}
\sqrt{n}\gamma_0 \P\bigl(L(0,0)> n\bigr) & =& \sum
_{k=0}^{\lfloor\sqrt
{n}\gamma_0 \rfloor} \E( \mathbf{1}_{\{L(k,n) > n\}} ) \times
\frac{\sqrt{n}\gamma_0}{\lfloor\sqrt{n}\gamma_0 \rfloor
+1}
\\
& =& \E( \xi_n) \times \frac{\sqrt{n}\gamma_0}{\lfloor\sqrt
{n}\gamma_0 \rfloor+1} \to
\E(\xi) = \frac{1}{\sqrt{\pi}}\qquad\mbox{as } n \to\infty.
\end{eqnarray*}
This proves Theorem \ref{clusterheight}.
\end{pf*}


Proposition \ref{propWeakConv1} will be proved through a sequence of lemmas.

To state the next lemma, we recall from Theorem \ref{theoremGRSDual-DobleBW}
that $(\wb{\mathcal X}_n,\wb{\wh{\mathcal X}}_n) \Rightarrow
({\mathcal W}, \widehat{{\mathcal W}})$ as $n \to\infty$.
Using Skorohod's representation theorem,
we assume that we are working on a probability space where
$d_{{\mathcal H}\times\widehat{{\mathcal H}}}((\wb{\mathcal X}_n,
\wb{\wh{\mathcal X}}_n), ({\mathcal W}, \widehat{{\mathcal
W}})) \to0 $ almost surely as $ n \to\infty$.
%
\begin{lemma}
\label{lemEtaptsetConv}
For $t_1 > t_0$, we have
\[
\P\bigl(\xi_{\wb{\mathcal X}_n} (t_0,t_1) \neq
\xi_{{\mathcal
W}}(t_0,t_1)\mbox{ for infinitely many }n \bigr) = 0.
\]
\end{lemma}

\begin{pf}
We prove the lemma for $t_0 = 0$ and $t_1 = 1$, that is, for $\xi_n =
\xi_{\wb{\mathcal X}_n} (0,1)$
and $\xi_{{\mathcal W}}(0,1)$, the proof for general $t_0,t_1$ being similar.
First, we show that, for all $k \geq0$,
%
\begin{equation}
\label{eqnweakconvliminf} \liminf_{n\to\infty} \mathbf{1}_{\{\xi_n \geq k\}} \geq
\mathbf {1}_{\{\xi\geq k\}}\qquad\mbox{almost surely}.
\end{equation}
Indeed, for $k = 0$, both $\mathbf{1}_{\{\xi_n \geq k\}}$ and
$\mathbf{1}_{\{\xi\geq k\}}$
equal $1$.
For $k\geq1$, (\ref{eqnweakconvliminf}) follows from
almost sure convergence of $(\wb{\mathcal X}_n,\wb{\widehat
{\mathcal X}}_n)$ to $
({\mathcal W}, \widehat{{\mathcal W}})$ and from the properties of the
set ${\mathcal M}_{\mathcal W}(0,1)$
as described in Proposition \ref{lemPropEtaPtset}.

To complete the proof, we need to show that $\P(\limsup_{n \to\infty
} \{\xi_n > \xi\}) = 0$.
This is equivalent to showing that $\P(\Omega^{k}_0) = 0$ for all
$k\geq0$,
where
\[
\Omega^{k}_0:= \bigl\{\omega: \xi_n(\omega)
> \xi(\omega) = k\mbox{ for infinitely many }n\bigr\}. %
\]
Consider $k=0$ first. From Proposition \ref{lemPropEtaPtset}, it
follows that on the event $\xi=0$,
almost surely we can obtain $\gamma:=\gamma(\omega)>0$ such that
${\mathcal M}_{{\mathcal W}}(0,1)
\cap(-\gamma, 1 + \gamma) = \varnothing$.
From the almost sure convergence of $(\wb{\mathcal X}_n,\wb
{\wh{\mathcal X}}_n)$ to $
({\mathcal W}, \widehat{{\mathcal W}})$,
we have $\P(\Omega^0_0) = 0$.

For $k > 0$, on the event $\Omega^{k}_0$ we show
a forward path $\pi\in{\mathcal W}$ coincides with a dual path
$\hat{\pi}
\in\widehat{\mathcal W}$ for a positive time which leads to a contradiction.
From Proposition \ref{lemPropEtaPtset}, it follows that given $\eta> 0$,
there exist $m_0 \in\N$ and $s_0 \in(1/m_0,1)$ such that $\P(\xi
_{\mathcal W}(1/m_0, 1)
= \xi_{\mathcal W}(1/m_0, s_0) = \xi_{\mathcal W}(0, 1) = k ) > 1-
\eta$,
that is, the paths leading to any single point considered in $
{\mathcal M}_{{\mathcal W}}(0,1) =
{\mathcal M}_{{\mathcal W}}(1/m_0,1)$ have coalesced before time $ s_0$.
Fix $ 0 < \varepsilon< 1/m_0 $ such that
$(x - \varepsilon,x + \varepsilon) \subset(0,1)$ for all $x \in{\mathcal
M}_{{\mathcal W}}(1/m_0, 1) $
and the $\varepsilon$-tubes around the $k$ paths contributing to
${\mathcal M}_{{\mathcal W}}(s_0, 1) $, {viz}., $\pi_1 (t),
\dotsc, \pi_k(t), t \in[s_0, 1]$, given by
\[
T^{i}_{\varepsilon}:= \bigl\{(x,t): \pi_i(t) -
\varepsilon\leq x \leq\pi _i(t) + \varepsilon, s_0 \leq t
\leq1 \bigr\}\qquad\mbox{ for }i=1,\dotsc,k, %
\]
are disjoint.
Since we have almost sure convergence on the event $\Omega^k_0$, there
exists $n_0$ such that
one of the $k$ tubes must contain at least two paths, $ \pi_1^{n_0},
\pi_2^{n_0}$ (say)
of ${\mathcal X}_{n_0}$ which do not coalesce by time $1$.
From the construction of dual paths,
it follows that there exists at least one dual
path $\hat{\pi}^{n_0} \in\wb{\widehat{{\mathcal X}}}^{1+}_{n_0}$
lying between $\pi_{1}^{n_0}$ and $\pi_{2}^{n_0}$ for $t\in[s_0,1]$,
and hence
we must have an approximating $\hat{\pi} \in\widehat{\mathcal
W}^{1+}$ close to
$\hat{\pi}^{n_0}$ for $t \in[s_0,1]$.
Since we have only finitely many disjoint $k$ tubes,
taking $\varepsilon\to0$ and using compactness of $\widehat{\mathcal
W}$ we obtain that there
exists $\hat{\pi} \in\widehat{\mathcal W}$ such that
$ \hat{\pi} (t) = \pi_i(t) $ for $ t \in[s_0, 1]$ and for some
$ 1\leq i \leq k$.
This violates the property of Brownian web and
its dual that they do not spend positive Lebesgue time together.
Hence, $\P(\Omega^{k}_{0}) = 0$ for all $k \geq0$ and this completes
the proof of the lemma.
\end{pf}

Lemma \ref{lemEtaptsetConv} immediately gives the following corollary.
%
\begin{cor}
\label{corEtaweakConv}
As $n \to\infty$, $\xi_n$ converges in distribution to $\xi$.
\end{cor}

Corollary \ref{corEtaweakConv} along with the following lemma
completes the proof of
Proposition \ref{propWeakConv1}.
%
\begin{lemma}
\label{lemUI}
The family $\{\xi_n: n \in\N\}$ is uniformly integrable.
\end{lemma}

\begin{pf}
For $m \in\N$, let
\[
K_m = [-m,m]^2 \cap{\mathbb Z}^2\quad
\mbox{and}\quad\Omega_m:= \bigl\{ (0,1),(0,-1),(1,1),(1,-1)\bigr
\}^{K_m}.
\]
We assign the product probability measure $\P^{\prime}$ whose
marginals for $ \bu\in K_m $ are
given by
\[
\P^{\prime}\bigl\{\zeta: \zeta(\bu) = (a,b)\bigr\} = \cases{
\displaystyle\frac{p}{2}, &\quad for $a = 1$ and $b \in\{1,-1\}$,
\vspace*{3pt}\cr
\displaystyle\frac{(1-p)}{2},
&\quad for $a = 0$ and $b \in\{1,-1\}$.}
\]
$\P^{\prime}$ is the measure induced by the random variables $\{
(B_{\bu}, U_{\bu}): \bu\in K_m \}$.

For ${\mathcal\zeta} \in\Omega_m$
and for $K \subseteq K_m $, 
the $K$ cylinder of $\zeta$ is given by $C(\zeta, K):=\{\zeta
^{\prime}: \zeta^{\prime}(\bu) =
\zeta(\bu)$ for all $\bu\in K\}$. For any two events $A, B \subseteq
\Omega_m$, let
\begin{eqnarray*}
A \Box B &:=& \bigl\{\zeta: \mbox{there exists }K = K(\zeta) \subseteq
K_m 
\mbox{ such that } C(\zeta, K) \subseteq A,
\\
&& \mbox{and }C\bigl(\zeta, K^{\prime}\bigr)\subseteq B\mbox{ for
} K^{\prime} = K_m 
\setminus K \bigr\}
\end{eqnarray*}
denote the disjoint occurrence of $A$ and $B$.
Note that this definition is associative, that is,
for any $A, B, C \subseteq\Omega_m$
we have $(A \Box B)\Box C = A \Box(B\Box C)$.


Let
\begin{eqnarray*}
F^m_n&:= & \bigl\{\mbox{there exist
}(u_1,n), (u_2, n) \in\wh{V} \mbox{ with } 0 \leq u_1 < u_2 \leq\sqrt{n}\gamma_0\mbox{ and}
\\
&&{} \bigl(v_{1}^l, l\bigr), \bigl(v_{2}^l,
l\bigr) \in V \mbox{ for all } 0 \leq l \leq n\mbox{ such that}
\\
&&{} {-}m \leq v_{1}^l < \hat{h}^l
(u_1,n) (1) < \hat{h}^l (u_2,n) (1) <
v_{2}^l \leq m \bigr\},
\\
E^m_n(k) &:=& \bigl\{\mbox{for } 1 \leq i \neq j\leq k,
\mbox{ there exists } (x_i,0) \in V \mbox{ with}
\\
&&{} h^n(x_i,0) (1) \in[0, \sqrt{n}\gamma_0 ]
, h^n(x_i,0) \neq h^n(x_j,0),
h^l(x_i, 0 ) (1) \in[-m,m]
\\
&&{} \mbox{for all }0\leq l \leq n\bigr\}.
\end{eqnarray*}
%

We claim that for all $k \geq2$,
%
\begin{equation}
\label{eqGRSBK} E^m_n(3k)\subseteq\underbrace{F^m_n
\square F^m_n\square \cdots\square F^m_n}_{k~\mathrm{times}}.
\end{equation}

We prove it for $k=2$. For general $k$, the proof is similar.
Let $(u_i, n) \in\wh{V}, 1 \leq i \leq5$ and $ (x_i, 0) \in V,
1 \leq i \leq6$ be
as in Figure~\ref{figUI}. The region explored to obtain the vertex $
\hat{h}^j ( u_i,n)$ for $ 1 \leq j \leq n $
is contained in $ \bigcup_{l=0}^{n-1} [h^l (x_{i},0)(1),\break h^l (x_{i+1},0)
(1) ]\times\{l\} $. Thus, the regions
explored to obtain the dual paths starting from $ (u_1,n), (u_2,n)$ and
the dual paths starting from $ (u_4,n), (u_5,n)$ are disjoint (see
Figure~\ref{figUI}).
Hence, it follows that $ E^m_n(6) \subseteq F^m_n \square F^m_n$.

%
\begin{figure}

\includegraphics{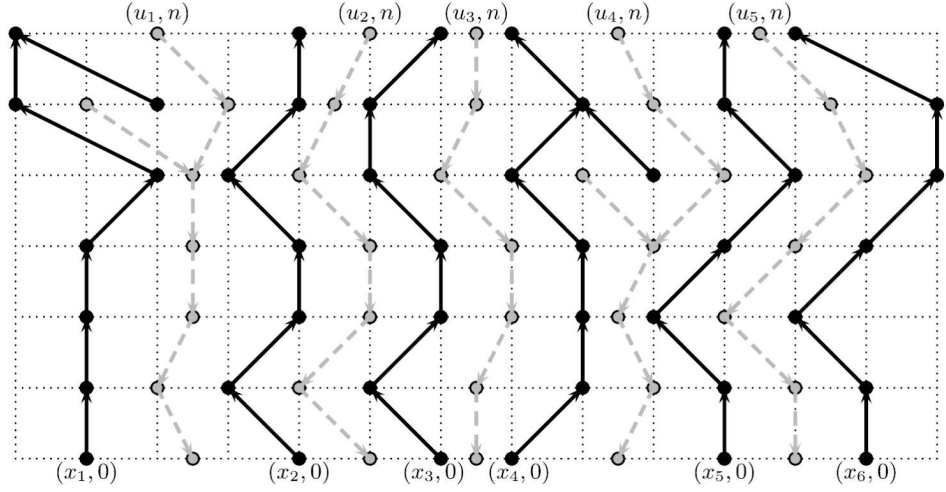}

\caption{The event $E^m_n(6)$.}
\label{figUI}
\end{figure}


Since the event $E^m_n(k)$ is monotonic in $m$, from (\ref{eqGRSBK})
we get
\begin{eqnarray*}
\P(\xi_n \geq3k) & =& \P\Bigl(\lim_{m \to\infty}
E^m_n(3k)\Bigr) = \lim_{m
\to\infty}\P\bigl(
E^m_n(3k)\bigr)
\\
& \leq& \lim_{m \to\infty}\P\bigl( F^m_n \Box
\cdots\Box F^m_n\bigr) = \lim_{m \to\infty}
\P^{\prime}\bigl( F^m_n \Box\cdots\Box
F^m_n\bigr).
\end{eqnarray*}
Applying the  BKR inequality (see Reimer \cite{R00}), we get
%
\begin{equation}
\label{eqUI} \P(\xi_n \geq3k) \leq\lim_{m \to\infty}
\bigl(\P^{\prime}\bigl( F^m_n \bigr)
\bigr)^k = \Bigl(\P\Bigl(\lim_{m \to\infty}
F^m_n \Bigr)\Bigr)^k = \bigl(
\P(F_n)\bigr)^k, 
\end{equation}
where $ F_n:= \{$there exist $ (u_1, n), (u_2, n) \in\wh{V} $ with
$ 0\leq u_1 < u_2 \leq\sqrt{n}\gamma_0 $ such that $
\hat{h}^{n} (u_1,n) \neq\hat{h}^{n} (u_2,n) \}$.


For any $(x,t) \in\R^2$ fix $ t_n
= \lfloor n t \rfloor$ and $ x_n = \max\{ \lfloor\sqrt{n} \gamma_0
x \rfloor+ j: j \leq0$,
$( \lfloor\sqrt{n} \gamma_0 x \rfloor+ j, t_n ) \in\wh{V} \}
$. Let
$ \hat{\theta}^{(x,t)}_n \in \widehat{{\mathcal X}}_n (\gamma
_0) $ be the scaling of the path
$\hat{\pi}^{(x_n,t_n)} \in\widehat{{\mathcal X}} $.
Define
\begin{eqnarray*}
F^{\prime}_n&:= & \bigl\{\hat{\theta}^{(0,1)}_n
\mbox{ and } \hat{\theta}^{(1,1)}_n \mbox{ do not
coalesce in time } 1\bigr\}.
\end{eqnarray*}

We observe that $ F_n \subseteq F^{\prime}_n $.
Now $\P(F^{\prime}_n)$ converges to the probability that two
independent Brownian
motions starting at a distance $1$ from each other do not meet
by time $1$. Since $\lim_{n\rightarrow\infty}\P(F^{\prime}_n)<1$,
the family $\{\xi_n:n\in\N\}$ is uniformly integrable.
\end{pf}

%
\begin{remark}
It is to be noted that Newman, Ravishankar and Sun \cite{NRS05} also used ideas of negative correlation
to establish the weak convergence of ${\mathcal M}_{\wb{\mathcal X}_n}$ as a point process on $\R$
for a more general setup where paths can cross each other.
In our case, the negative correlation ideas come in a much less
essential manner only
to establish uniform integrability as the noncrossing nature of paths
enable us to
obtain Corollary \ref{corEtaweakConv}. 
\end{remark}


\section{Proofs of Theorems \texorpdfstring{\protect\ref{BM}}{1.3} and \texorpdfstring{\protect\ref{BEarea}}{1.4}}\label{BMBEarea}

%
\begin{figure}[b]

\includegraphics{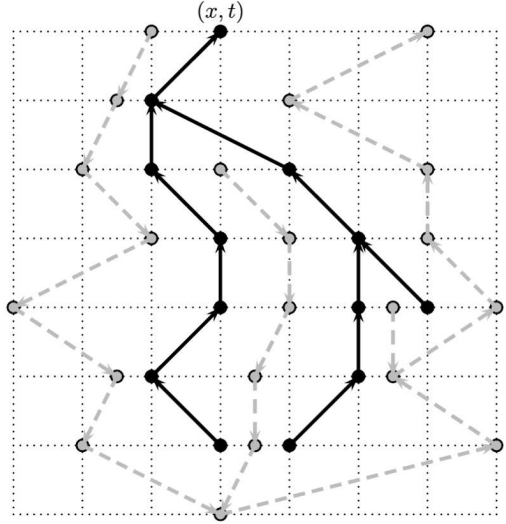}

\caption{The two dual paths $\hat{\pi}^{\hat{l}(x,t)}$ and
$\hat{\pi}^{\hat{r}(x,t)}$ enclose the cluster $C(x,t)$.
These dual paths after scaling are each Brownian paths.}
\label{ideafigure}
\end{figure}


In this section, we prove Theorems \ref{BM} and \ref{BEarea}.
The main idea of the proof is that the horizontal distance between the
dual paths $\hat{\pi}^{\hat{r}(x,t)}$ and $\hat{\pi
}^{\hat{l}(x,t)}$ (see Figure~\ref{ideafigure}) form a Brownian
excursion process after scaling. The cluster $C(x,t)$ being enclosed
between these two paths, its size is related to the area under the
Brownian excursion.

For a formal proof, we need to introduce some notation.
For $\tau> 0$, let $S^\tau,S^{\tau^+}:C[0,\infty) \to\R$ be
defined by
$S^\tau(f):= \inf\{t \geq0: f(t + s)\geq f(t )$ for all $0\leq s
\leq\tau\}$
and $S^{\tau^+}(f):= \inf\{t \geq0: f(t + s) > f(t )$ for all $0 <
s \leq\tau\}$.
Let $T^{\tau^+}: C[0,\infty) \to C[0,\infty)$ be the map given by
%
\begin{eqnarray}
\label{eqfMeander} T^{\tau^+}(f) (s):= \cases{ f\bigl(S^{\tau^+} + s
\bigr)-f\bigl(S^{\tau^+}\bigr), &\quad if $S^{\tau^+} < \infty$,
\cr
f(s), &
\quad otherwise.}
\end{eqnarray}
For a Brownian motion $ W $ with $ W(0) = 0$, we define $W^{\tau} =
T^{\tau^+}(W)$.
From Bolthausen \cite{B76}, we have $S^{\tau^+} = S^\tau< \infty$ almost surely
under the measure induced by $W$ on $C[0,\infty)$ and
$W^{1}\mid _{[0,1]} \disteq W^{+}$ where $W^+$ is the standard Brownian
meander process defined in (\ref{eqBMBE}). From the scaling
property of Brownian motion, it follows that
$ \{ W^{\tau} (s): s \in[0, \tau] \}
\disteq\{ \sqrt{\tau}W^{+}(s / \tau): s \in[0, \tau] \} $.
Durrett, Iglehart and Miller \cite{DIM77} (Theorem 2.1) proved that
$ W \mid {\mathbf1}_{ \{\min_{s \in[0,1]}W (s) \geq-\varepsilon\} } \weak
W^{+}$ as $ \varepsilon\downarrow0$. Using this result and the scaling
property of $W^{\tau}$, given above,
straightforward calculations imply the following lemma and its
corollary (for details, see Roy, Saha and Sarkar \cite{RSS15}).
%
\begin{lemma}
\label{lemRandomMeander}
For $\tau> 0$ considering $W$ as a standard Brownian motion on
$[0,\infty)$ starting from $0$, we have
$W \mid {\mathbf1}_{\{\min_{t\in[0,\tau]}W(t) \geq-1/n\} }\weak
W^{\tau}$ as $n \to\infty$.
\end{lemma}
Define $\widetilde{W}^{\tau}$ as the process on $C[0,\infty)$ given by
\begin{eqnarray*}
\widetilde{W}^{\tau}(t):= \cases{ W^{\tau}(t), &\quad if $0 \leq t
\leq\tau$,
\cr
W^{\tau}(\tau) + \widetilde{W}(t-\tau), &\quad otherwise,}
\end{eqnarray*}
where $\widetilde{W}$ is a Brownian motion on $[0,\infty)$,
independent of $W^{\tau}$, with $\widetilde{W} (0) = 0$.
For $f \in C[0,\infty)$, let $t_f:= \inf\{s > 0: f(s) = 0\}$
with $t_f = \infty$ if $f(s)\neq0$ for all $s > 0$.
Consider the mapping $H:C[0,\infty) \to C[0,\infty)$ given by $
H(f)(t):= \mathbf{1}_{\{t \leq t_f\}}f(t)$.
We define $W^{+,\tau} = H(W^{\tau})$.
A
similar argument as that of Lemma \ref{lemRandomMeander} gives us
the following corollary.
%
\begin{cor}
\label{corRandomBetaMeander}
For $\tau> 0$, we have, $W^{\tau} \disteq\widetilde{W}^{\tau}$ and
$W^{+,\tau}
\disteq H(\widetilde{W}^{\tau})$.
\end{cor}
Let $A\subset C[0,\infty)$ be such that
%
\begin{eqnarray}
\label{defAset} A &:= & \bigl\{f \in C[0,\infty): t_f <\infty\mbox{
and for every } \varepsilon> 0 \mbox{ there exists }
\nonumber\\[-8pt]\\[-8pt]\nonumber
&&{} s\in(t_f, t_f+\varepsilon) \mbox{ with }f(s)<0\bigr\}.
\end{eqnarray}
From Corollary \ref{corRandomBetaMeander}, it follows that
$\P(W^{\tau} \in A) = 1$. Hence, $H$ is continuous almost surely
under the measure induced by $W^\tau$ on $C[0,\infty)$.

Next, we obtain the distribution of $\int_0^\infty W^{+,\tau}(t)\,dt$.
%
\begin{lemma}
\label{lemmaRandomExcursionArea}
For $\tau, \lambda > 0$, we have
\begin{eqnarray*}
&& \P\biggl(\int_0^\infty W^{+,\tau}(t)\,dt >
\lambda\biggr) = \frac{ \sqrt{\tau}}{2} \int_\tau^{\infty}
t^{-3/2}\wb{F}_{I^+_0} \bigl(\lambda t^{-3/2}\bigr)\,dt.
\end{eqnarray*}
\end{lemma}

\begin{pf}
We give here a straightforward
proof using random walk.
Let $\{S_n: n \geq0\}$ be a
symmetric random walk with variance $1$ starting at $S_0 = 0$.
Since $\P(W^\tau\in A) = 1$, minor modification of the argument used
to prove Lemmas 2.4 and
2.5, Bolthausen \cite{B76} shows that $H\circ T^{\tau^+}$ is almost surely
continuous under
the measure induced by $W$ on $C[0,\infty)$ (for details, see Roy, Saha and Sarkar \cite{RSS15}).
From Donsker's invariance principle and from the
continuous mapping theorem, it follows that
for $\lambda> 0$, a continuity point of
$\int_0^\infty W^{+,\tau}(t) \,dt$, we have
%
\begin{eqnarray*}
&&\P\biggl(\int_0^\infty W^{+,\tau}(t)\,dt >
\lambda\biggr) = \lim_{n \to\infty} \P\biggl(\int_0^\infty
H\bigl(T^{\tau^+}(Y_n)\bigr) (t)\,dt > \lambda\biggr),
\end{eqnarray*}
where
%
\begin{equation}
\label{eqYn} Y_n (t):= \frac{ S_k }{\sqrt{n} } + \frac{ (nt - [nt]) }{ \sqrt{n}} (
S_{k+1} - S_k )\qquad\mbox{for }\frac{ k}{ n} \leq t <
\frac{ k+1}{n}.
\end{equation}
A
similar argument as in Lemma 3.1 of  Bolthausen \cite{B76} gives us that
(for details, see Roy, Saha and Sarkar \cite{RSS15})
\begin{eqnarray*}
&& \P\biggl(\int_0^\infty H\bigl(T^{\tau^+}(Y_n)
\bigr) (t)\,dt > \lambda\biggr)
\\
&&\qquad = \P\biggl(\int_0^\infty H(Y_n)
(t)\,dt > \lambda\Big| \min_{t\in[0,\tau
]}Y_n(t) \geq0,
t_0 > n\tau\biggr),
\end{eqnarray*}
where
$t_0:= \inf\{n>0: S_n =0\}$ is the first return time to $0$ of the
random walk.
Hence for $\lambda> 0$, a continuity point of $W^{+,\tau}$, we
obtain 
%
\begin{eqnarray*}
&& \P\biggl(\int_0^\infty W^{+,\tau}(t)\,dt >
\lambda\biggr)
\\
&&\qquad  = \lim_{n \to\infty} \P\biggl(\int_0^\infty
H(Y_n) (t)\,dt > \lambda \Big| \min_{t\in[0,\tau]}Y_n(t)
\geq0, t_0 > n\tau\biggr)
\\
&&\qquad  = \lim_{n \to\infty} \sum_{j=1}^{\infty}
\frac{n^{3/2}\P(t_0 = \lfloor n\tau\rfloor+ j)}{n(\sqrt
{n}\P(t_0 > n\tau))}
\\
&&\quad\qquad{} \times\P\biggl(\int_0^\infty
H(Y_n) (t)\,dt > \lambda \Big| \min_{t\in[0,\tau]}Y_n(t)
\geq0, t_0 = \lfloor n\tau\rfloor+ j\biggr)
\\
&&\qquad  = \lim_{n \to\infty} \frac{1}{\sqrt{n}\P(t_0 > n \tau)}\int_{\lfloor n\tau
\rfloor/n}^\infty
g_n(t) f_n(t)\,dt,
\end{eqnarray*}
where\vspace*{1pt} for $t\geq\lfloor n\tau\rfloor/n$,
$f_n(t) = \P(\int_0^\infty H(Y_n)(u)\,du > \lambda\mid  \min_{t\in
[0,\tau]}Y_n(t) \geq0,
t_0 = \lfloor n t \rfloor+ 1)$ and
$g_n(t) = n^{3/2}\P(t_0 = \lfloor n t \rfloor+ 1)$.
It is known that (see Kaigh \cite{K76})
\[
\lim_{n\to\infty}\sqrt{n}\P(t_0 > n) = \sqrt{
\frac{2}{\pi
}}\quad\mbox{and}\quad \lim_{n\to\infty}n^{3/2}
\P(t_0 = n) = \frac{1}{\sqrt
{2\pi}}. %
\]
Hence, from Theorem 2.6 Kaigh \cite{K76} together with
the continuous mapping theorem and the
scaling property of the Brownian motion we have
$\P(\int_0^\infty W^{+,\tau}(t)\,dt > \lambda) = \frac{ \sqrt{\tau
}}{2} \int_\tau^{\infty} t^{-3/2}\wb{F}_{I^+_0}(\lambda
t^{-3/2})\,dt$.
Finally, $I^+_0$ being a continuous random variable (see Louchard and Janson \cite
{JL07}), it
follows that the random variable $\int_0^\infty W^{+,\tau}(t)\,dt$ is
continuous. This
completes the proof.
\end{pf}

\subsection{Proof of Theorem \texorpdfstring{\protect\ref{BM}}{1.3}}
Recall that $\hat{r}(x,t)$ and $\hat{l}(x,t)$ denote
the right and left dual neighbours, respectively, of $(x,t)\in V$. Let
$\hat{D}_k(x,t):=
\hat{h}^k (\hat{r}(x,t))(1) - \hat{h}^k (\hat{l}(x,t))(1)$ where $\hat{h}$ is as defined after (\ref
{eqndefalar}). Consider the continuous
function $\hat{D}^{(x,t)}_n \in C[0,\infty)$ given by
%
\begin{eqnarray}\label{eqDualDn}
&&\hat{D}_n^{(x,t)} (s):=
\frac{ \hat{D}_k (x,t) }{ \gamma_0
\sqrt{n} } + \frac{ (ns - [ns]) }{ \gamma_0 \sqrt{n} } \bigl( \hat{D}_{k+1} (x,t) -
\hat{D}_k (x,t)\bigr)
\nonumber\\[-8pt]\\[-8pt]
\eqntext{\displaystyle\mbox{for }\frac{ k}{ n} \leq s
\leq\frac{ k+1}{n}.}
\end{eqnarray}

Fix $\tau> 0$.
For an ${\mathcal H}\times\widehat{{\mathcal H}}$ valued random
variable $(K,\widehat{K})$ and
for $x \in{\mathcal M}_{K}(0, \tau)$ let $\hat{\pi}^{(x,\tau
)}_r$ be defined as
\begin{eqnarray*}
\hat{\pi}^{(x,\tau)}_r:= \cases{ \hat{\pi}, &\quad if $
\sigma_{\hat{\pi}} = \tau$ and there is no $\hat{\pi}_1 \in
\widehat{K}^{\tau+}$ with $x < \hat{\pi}_1(\tau)<\hat{\pi}(\tau)$,
\cr
\hat{\pi}_0, &\quad otherwise,}
\end{eqnarray*}
where $\hat{\pi}_0$ denotes the constant zero function with
$\sigma_{\hat{\pi}_0} = \tau$. In other words,
$\hat{\pi}^{(x,\tau)}_r \in\widehat{K}^{\tau+}$ is such that
among all $\hat{\pi}\in\widehat{K}^{\tau+}$,
$\hat{\pi}^{(x,\tau)}_r(\tau)$ is closest to $(x,\tau)$
on the right. Similarly, $\hat{\pi}^{(x,\tau)}_l$ is defined as
the path closest to $(x,\tau)$
on the left.\vspace*{2pt}

For $\hat{\pi}\in\wh{\Pi}$ with $\sigma_{\hat{\pi}}
\geq\tau$, let $g(\hat{\pi}) \in C[0,\infty)$ be given by $
g(\hat{\pi})(t):= \hat{\pi}(\tau-t)$ for $t \geq0$.
Fix $f \in C_b[0,\infty)$ and define
\begin{eqnarray*}
&&\kappa_{ (K, \widehat{K}) }(\tau,f):= \sum_{x \in{\mathcal M}_{K}(0,\tau)} f
\bigl(g\bigl(\hat{\pi}^{(x,\tau
)}_r\bigr) - g\bigl(\hat{\pi}^{(x,\tau)}_l\bigr)\bigr).
\end{eqnarray*}
Let\vspace*{1pt} $\kappa(\tau,f):= \kappa_{({\mathcal W},\widehat{\mathcal
W})}(\tau,f)$,
and $\kappa_n(\tau,f):= \kappa_{(\wb{\mathcal X}_n,\wb{\widehat
{{\mathcal X}}}_n)}(\tau,f) $.
Comparing with the definitions introduced
in (\ref{eqnDefNK}), for $m_f = \sup\{\llvert   f(s)\rrvert: s \in[0,\infty)\}$
we have
%
\begin{equation}
\label{eqnRelationKappa} \kappa(\tau,f) \leq m_f\xi_{{\mathcal W}}(0,\tau),
\kappa_n(\tau,f) \leq m_f\xi_{\wb{\mathcal X}_n} (0,\tau)
\qquad\mbox{for all } n \geq1.
\end{equation}
From Proposition \ref{lemPropEtaPtset}, we know that for each
$x \in{\mathcal M}_{\mathcal W}(0,\tau)$, there exist $\hat{\pi
}^{(x, \tau)}_r$,
$\hat{\pi}^{(x, \tau)}_l \in\widehat{\mathcal W}$ both starting
from $ (x, \tau)$
with $\hat{\pi}^{(x, \tau)}_r(0) >
\hat{\pi}^{(x, \tau)}_l(0) $.

The following lemma is the main tool for
establishing Theorem \ref{BM} and Theorem \ref{BEarea}.
%
\begin{lemma}
\label{lemmaKappanKappaExp}
For $\tau> 0$ and $f \in C_b[0,\infty)$, we have
%
\begin{equation}
\label{eqnKnLimitK} \lim_{n \to\infty}\E\bigl[\kappa_n(\tau,f)
\bigr] = \E\bigl[\kappa(\tau,f)\bigr].
\end{equation}
\end{lemma}

\begin{pf} 
From (\ref{eqnRelationKappa}) and Lemma \ref{lemUI}, it follows that
the family $\{\kappa_n(\tau,f): n \in\N\}$ is uniformly integrable.
Hence, it suffices to show that $\kappa_n(\tau,f)$ converges in
distribution to $\kappa(\tau,f)$ as
$n \to\infty$. We assume\vspace*{1pt} that we are working on
a probability space such that $(\wb{{\mathcal X}}_n, \wb{\widehat
{{\mathcal X}}}_n)$ converges
to $({\mathcal W}, \widehat{{\mathcal W}})$ almost surely in
$({\mathcal H}\times\widehat{\mathcal H},
d_{{\mathcal H}\times\widehat{{\mathcal H}}})$.
From Lemma \ref{lemEtaptsetConv},
we have $\lim_{n \to\infty}\xi_{\wb{{\mathcal X}}_n}(0,\tau) =
\xi_{\mathcal W}(0,\tau)$ almost surely,
and hence from (\ref{eqnRelationKappa}) for $ \xi_{\mathcal
W}(0,\tau) = 0$, we have $
\kappa_n(\tau,f) = \kappa(\tau,f) = 0$ for all $ n $ large.
Next, we consider the case $\xi_{\mathcal W}(0,\tau) = k\geq1$.
Suppose ${\mathcal M}_{\mathcal W}(0,\tau) = \{x_1,\ldots, x_k \}$.
From Lemma \ref{lemEtaptsetConv},
we have\vspace*{2pt} that ${\mathcal M}_{\wb{{\mathcal X}}_n}(0,\tau) = \{
x^n_1,\ldots, x^n_k \}$
for all large $n$ and $\lim_{n \to\infty} x^n_i = x_i$ for all $1
\leq i \leq k$.
Fix $T\geq0$. To complete the proof, it is enough to show that
$\sup\{\llvert   \hat{\pi}^{(x_i,\tau)}_r(\tau-s) - \hat{\pi
}^{(x^n_i,\tau)}_r(\tau-s)\rrvert
\vee\llvert   \hat{\pi}^{(x_i,\tau)}_l(\tau-s) - \hat{\pi
}^{(x^n_i,\tau)}_l(\tau- s)\rrvert:
s \in[0, \tau+ T]\} \to0$ as $n \to\infty$ for all $1 \leq i \leq k$.

We observe that for $y_i \in(\hat{\pi}^{(x_i,\tau)}_r(0),
\hat{\pi}^{(x_i,\tau)}_l(0))
\cap\Q$ there exists $\pi^{(y_i,0)} \in{\mathcal W}$ such that $\pi
^{(y_i,0)}(\tau) = x_i$.
We choose $ \varepsilon= \varepsilon(\omega) > 0 $ so that for all $1 \leq
i \leq k$:
\begin{longlist}[(a)]
\item[(a)] $(x_i- \varepsilon, x_i +\varepsilon) \subset(0,1)$, $(x_i-
2\varepsilon, x_i +2\varepsilon)
\cap{\mathcal M}_{\mathcal W}(0, \tau) = \{x_i\}$ and
\item[(b)] $(\hat{\pi}^{(x_i,\tau)}_r(0)-\pi
^{(y_i,0)}(0))\wedge
(\pi^{(y_i,0)}(0)-\hat{\pi}^{(x_i,\tau)}_l(0)) > 2 \varepsilon$.
\end{longlist}
Let $n_0 = n_0(\omega) $ be such that,
for all $n \geq n_0$:
\begin{longlist}[(ii)]
\item[(i)] $\xi_{\wb{{\mathcal X}}_n}(0,\tau) = \xi_{\mathcal
W}(0,\tau)$
and
\item[(ii)] for all $1 \leq i \leq k$ there exist $\hat{\pi
}^{1,n}_i, \hat{\pi}^{2,n}_i
\in\wb{\wh{\mathcal X}}^{\tau+}_n$ and
$\pi^{n}_i \in\wb{{\mathcal X}}^{0 -}_n$ such that
$\sup\{\llvert   \hat{\pi}^{1,n}_i(\tau- s) - \hat{\pi}^{(x_i,\tau
)}_r(\tau- s)\rrvert    \vee
\llvert   \hat{\pi}^{2,n}_i(\tau- s) - \hat{\pi}^{(x_i,\tau
)}_l(\tau- s)\rrvert    \vee
\llvert   \pi^{n}_i(\tau- s) - \pi^{(y_i,0)}(\tau- s)\rrvert: s \in[0,\tau+
T]\} < \varepsilon$.
\end{longlist}
The choice of $n_0 $ ensures that ${\mathcal M}_{\wb{{\mathcal X}}_n}(0,\tau) \cap(x_i - \varepsilon,
x_i + \varepsilon) = \{x^{n}_i\}$. Since there exist only two dual paths
starting from $(x_i,\tau)$,
because of the uniqueness of $x_i^n$ in the interval
$(x_i-\varepsilon, x_i+\varepsilon)$ and the noncrossing nature of our
paths we must have
$\hat{\pi}^{(x^n_i,\tau)}_r(\tau-s) = \hat{\pi
}^{1,n}_i(\tau-s)$ and
$\hat{\pi}^{(x^n_i,\tau)}_l(\tau-s) = \hat{\pi
}^{2,n}_i(\tau-s)$ for all
$s \in[0, \tau+T]$ and for all $n \geq n_0$ (for details, see Roy, Saha and Sarkar \cite{RSS15}).
Since $T\geq0$ is chosen arbitrarily, this completes the proof.
\end{pf}

The next lemma calculates $\E[\kappa(\tau,f)]$.
%
\begin{lemma}
\label{lemmaOpenSetBdry}
For $\tau>0$ and $f \in C_b[0,\infty)$, we have
\[
\E\bigl[\kappa(\tau,f)\bigr] = \E\bigl( f\bigl(\sqrt{2} W^{+,\tau} \bigr)
\bigr) /\sqrt{\pi \tau}. %
\]
\end{lemma}

\begin{pf} Let\vspace*{2pt} $I_n \subset\{0,1,\ldots,n-1\}$ given by
$I_n:=
\{i: 0 \leq i \leq n-1, \hat{\pi}^{(i/n,\tau)},\break \hat{\pi
}^{((i+1)/n,\tau)} \in
\widehat{\mathcal W}$ such that $\hat{\pi}^{(i/n,\tau)}(0) <
\hat{\pi}^{((i+1)/n,\tau)}(0)\}$. We define
\[
{\mathcal R}_n(\tau,f) = \sum_{i \in I_n} f
\bigl(g\bigl(\hat{\pi }^{((i+1)/n,\tau)} - \hat{\pi}^{(i/n,\tau)}\bigr)
\bigr).
\]
From Proposition \ref{lemPropEtaPtset}, we know ${\mathcal
M}_{{\mathcal W}}(0,\tau) \cap\Q= \varnothing$.
For each $x \in{\mathcal M}_{\mathcal W}(0,\tau)$, set $ l^x_n
=\lfloor nx \rfloor/n $
and $ r^x_n = l^x_n + (1/n)$.
Since there are exactly two dual paths $\hat{\pi}^{(x,\tau)}_r$
and $\hat{\pi}^{(x,\tau)}_l$
starting from $(x,\tau)$ with $\hat{\pi}^{(x,\tau)}_r(0) >
\hat{\pi}^{(x,\tau)}_l(0)$,
from Proposition 3.2(e) of Sun and Swart \cite{SS08} it follows that
$\{\hat{\pi}^{(l^x_n,\tau)}: n \in\N\}$ and $\{\hat{\pi
}^{(r^x_n,\tau)}: n \in\N\}$
converge to $\hat{\pi}^{(x,\tau)}_l$ and $\hat{\pi
}^{(x,\tau)}_r$, respectively,
in $(\wh{\Pi}, d_{\wh{\Pi}})$ as $n \to\infty$.
Hence, $ {\mathcal R}_n(\tau,f) \to\kappa(\tau,f) $ almost surely\vspace*{1pt}
as $n \to\infty$.
For each $i \in I_n$, there exist $y_i \in(\hat{\pi}^{(i/n,\tau)}(0),
\hat{\pi}^{((i+1)/n,\tau)}(0))\cap\Q$ and $\pi^{(y_i,0)}\in
{\mathcal W}$ such that
$\pi^{(y_i,0)}(\tau) \in{\mathcal M}_{\mathcal W}(0,\tau)$.
Hence, for $m_f = \sup\{\llvert   f(t)\rrvert: t \geq0\}$ we have
${\mathcal R}_n(\tau,f) \leq m_f\xi_{\mathcal W}(0,\tau)$ for all
$n$. As
$\E[\xi_{\mathcal W}(0,\tau)] < \infty$, the family $\{{\mathcal
R}_n(\tau,f): n \in\N\}$
is uniformly integrable, and hence we have $\lim_{n \to\infty} \E
[{\mathcal R}_n(\tau,f)] =
\E[\kappa(\tau,f)]$. From the fact that $ g( \hat{\pi
}^{((i+1)/n,\tau)})
- g( \hat{\pi}^{(i/n,\tau)}) \disteq H(1/n + \sqrt{2}W) $
where $ W$ denotes the standard Brownian motion on $[0,\infty)$, we have
\begin{eqnarray*}
&& \lim_{n \to\infty} \E\bigl[{\mathcal R}_n(\tau,f)\bigr]
\\
&&\qquad = \lim_{n \to\infty}n \E\Bigl[f\bigl(H(1/n + \sqrt{2}W)\bigr) \big|
1/n + \min_{t\in[0,\tau]}\sqrt{2}W (t) > 0\Bigr]
\\
&&\quad\qquad{} \times\P\Bigl(1/n + \min_{t\in[0,\tau]}\sqrt{2}W (t) > 0
\Bigr)
\\
&&\qquad  = \lim_{n \to\infty}\E\Bigl[ f\bigl(H(1/n+ \sqrt{2}W) \bigr) \big|
\min_{t\in[0,\tau]}\sqrt{2}W(t) > -1/n \Bigr] n
\\
&&{}\quad\qquad{}\times \bigl(2\Phi(1/\sqrt{2
\tau }n) - 1\bigr)
\nonumber
\\
&&\qquad  = \E\bigl( f\bigl(\sqrt{2} W^{+,\tau}\bigr) \bigr) /\sqrt{\pi\tau},
\end{eqnarray*}
where the last equality follows from Lemma \ref{lemRandomMeander},
Slutsky's theorem and
continuous mapping theorem.
This completes the proof.
\end{pf}

Now, to complete the proof of Theorem \ref{BM} we need the following lemmas.
%
\begin{lemma}
\label{lemDhatBM}
For $\tau> 0$, we have
$\hat{D}^{(0,0)}_n \mid  {\mathbf1}_{\{ L(0,0) > n\tau\}}
\weak \sqrt{2} W^{+,\tau}$
as $n \to\infty$.
\end{lemma}

\begin{pf}
Using translation invariance of our model, we have
\[
\E\bigl(f\bigl(\hat{D}^{(0,0)}_n\bigr) \mid {
\mathbf1}_{\{ L(0,0) > n\tau\}}\bigr) = \frac{ \E[\kappa_n(\tau,f)] }{ \E[\xi_{\wb{\mathcal X}_n}(0,\tau)] } \to\frac{ \E[\kappa(\tau,f)] } { \E[\xi_{\mathcal W}(0,\tau)] } = \E
\bigl( f\bigl(\sqrt{2} W^{+,\tau} \bigr)\bigr). %
\]
This holds for all $f \in C_b[0, \infty)$ which completes the
proof.
\end{pf}

%
\begin{lemma}
\label{lemDhatDnKn}
For $\tau> 0$, we have:
\begin{longlist}[(a)]
\item[(a)] $\sup\{\llvert   \hat{D}^{(0,0)}_n(s) -
D^{(0,0)}_n(s)\rrvert: s \geq0\}
\mid  {\mathbf1}_{\{ L(0,0) > n\tau\}} \prob0$ as $n \to\infty$,

\item[(b)] $\sup\{\llvert   K^{(0,0)}_n(s) - pD^{(0,0)}_n(s)\rrvert: s \geq
0\}
\mid  {\mathbf1}_{\{ L(0,0) > n\tau\}} \prob0$ as $n \to\infty$.
\end{longlist}
\end{lemma}
%

\begin{pf}
For part (a), fix $0 < \alpha< 1/2$, $T
\geq0$ and we observe that
\begin{eqnarray*}
&& \P\bigl(\sup\bigl\{\bigl\llvert \hat{D}_k(0,0)-
D_k(0,0)\bigr\rrvert: k \geq0\bigr\} \geq n^{\alpha
}, L(0,0)
> n\tau\bigr)
\\
&&\qquad \leq\P\bigl(\max\bigl\{\bigl\llvert \hat{D}_k(0,0)-
D_k(0,0)\bigr\rrvert: 0 \leq k \leq n(\tau + T) + 1\bigr\} \geq
n^{\alpha},
\\
&&\quad\qquad{}  L(0,0) > n\tau\bigr)
+ \P\bigl(L(0,0) > n(\tau+T)\bigr).
\end{eqnarray*}
Because of Theorem \ref{clusterheight}, it is enough to show that
$ \sqrt{n}  \P(\max\{\llvert   \hat{D}_k(0,0)- D_k(0,0)\rrvert:0\leq k \leq
n(\tau+ T)+1\} \geq
n^{\alpha}, L(0,0) > n\tau)\to0 $ as $ n\to\infty$.
Here, we present the simple idea behind the proof; the details are
available in Roy, Saha and Sarkar \cite{RSS15}.

The distance $d^l_k$ between $l_k(0,0)$ and the closest open vertex to
the left of $l_k(0,0)$ being $n^{\alpha}$ or more has a probability
$(1-p)^{n^{\alpha}}$. Thus, the probability that the maximum such
difference for $0\leq k \leq n(\tau+ T)+1$
is bigger that $n^{\alpha}$ is of the order $n(1-p)^{n^{\alpha}}$.
Similarly, for the distance $d^r_k$ associated with the vertex $r_k(0,0)$.
Since $\llvert   \hat{D}_k(0,0)- D_k(0,0)\rrvert    \leq d^l_k + d^r_k$, as $ n\to
\infty$, we have that
$ \sqrt{n}  \P(\max\{\llvert   \hat{D}_k(0,0)- D_k(0,0)\rrvert:0\leq k \leq
n(\tau+ T)+1\} \geq
n^{\alpha}, L(0,0) > n\tau)$ converges to $0$.

For part (b) of the lemma, we need ${D}^{(0,0)}_n \mid  {\mathbf
1}_{\{ L(0,0) > n\tau\}} \weak \sqrt{2} W^{+,\tau}$
as $n \to\infty$ which follows from part (a) and Lemma \ref
{lemDhatBM}. Hence, $r_k(0,0) - l_k(0,0)$ is of the order $\sqrt{n}$.
Also given $l_k(0,0)$ and $r_k(0,0)$, the number of open vertices lying
between these vertices has a binomial distribution with parameters
$(r_k(0,0) - l_k(0,0) -1)$ and $p$. Since these open vertices together
with $l_k(0,0)$ and $r_k(0,0)$ constitute $C_k(0,0)$, the proof follows
from similar order comparisons as done in (a). 
\end{pf}

\begin{pf*}{Proof of Theorem \ref{BM}}
We remarked that $W^{1}\mid _{[0,1]} = W^{+,1}\mid _{[0,1]} \disteq W^{+}$.
The proof of Theorem \ref{BM} follows from Lemmas \ref{lemDhatBM}
and \ref{lemDhatDnKn} and Slutsky's theorem
with the choice of $\tau= 1$.
\end{pf*}

\subsection{Proof of Theorem \texorpdfstring{\protect\ref{BEarea}}{1.4}}

For $\lambda> 0$, let $\bar{\lambda}:= \lambda^{3/2}(\sqrt{2}\gamma_0p)^{-1}$.
We show that:
%
\begin{lemma}
\label{lemmaLambdaBetaPositiveNonTruncate}
For $ \tau, \lambda > 0 $,
\begin{eqnarray*}
&& \lim_{n \to\infty}\sqrt{n}\P \Biggl( L(0,0) > n\tau, \sum
_{k = 0}^
{\infty}\# C_k(0,0) > (\lambda
n)^{3/2} \Biggr)
\\
&&\qquad  = \frac{1}{ \gamma_0 \sqrt{\pi\tau}} \P \biggl( \sqrt{2}\int_{0}^{\infty}
W^{+,\tau}(t)\,dt > \bar{\lambda} \biggr)
\\
&&\qquad = \frac{1}{ 2\gamma_0 \sqrt{\pi}}\int
_\tau ^\infty\wb{F}_{I^+_0} \bigl(\bar{\lambda} t^{-3/2}\bigr)t^{-3/2} \,dt.
\end{eqnarray*}
\end{lemma}

\begin{pf} For $f \in C[0,\infty)$ let
$I(f):= \int_0^\infty H(f)(t)\,dt$. Since $\P(W^\tau\in A) = 1$ where
$A$ is defined
as in (\ref{defAset}), $I$ is almost surely continuous under the
measure induced by $W^\tau$ on
$C[0,\infty)$. The proof follows from
Theorem \ref{BM}(ii) and the continuous mapping theorem.
\end{pf}


From the previous lemma,
we derive the following.
%
\begin{cor}
\label{corLambdaPositiveNonTruncate}
For $ \lambda> 0 $, we have
\begin{eqnarray*}
&& \lim_{n \to\infty}\sqrt{n}\P \bigl(\# C(0,0) > (\lambda
n)^{3/2} \bigr) = \frac{1}{ 2\gamma_0 \sqrt{\pi}}\int_0^\infty
\wb{F}_{I^+_0} \bigl(\bar{\lambda} t^{-3/2}\bigr)t^{-3/2}
\,dt.
\end{eqnarray*}
\end{cor}

\begin{pf} For any $\tau> 0$, we have
$\P(\# C(0,0) > (n\lambda)^{3/2}) \geq\P(L(0,0)>n\tau,\# C(0,0) >
(n\lambda)^{3/2})$, and hence\vspace*{2pt}
$\liminf_{n \to\infty}\sqrt{n}\P(\# C(0,0) > (n\lambda)^{3/2})
\geq
\frac{1}{ 2\gamma_0 \sqrt{\pi}}\int_0^\infty\wb{F}_{I^+_0}
(\bar{\lambda} t^{-3/2})t^{-3/2} \,dt$.\vspace*{2pt}

We observe that
\begin{eqnarray*}
&& \sqrt{n}\P\bigl(L(0,0)\leq n\tau, \# C(0,0) > (n\lambda)^{3/2}\bigr)
\\
&&\qquad  \leq\sqrt{n}\P\Biggl(\sum_{k=0}^{\lfloor n\tau\rfloor}
\hat{D}_k(0,0) > (n\lambda)^{3/2}\Biggr)
\\
&&\qquad  \leq\sqrt{n}\E\Biggl[\sum_{k=0}^{\lfloor n\tau\rfloor}
\widehat {D}_k(0,0)\Biggr](n\lambda)^{- 3/2}
\\
&&\qquad  = \sqrt{n}\bigl(\lfloor n\tau\rfloor+ 1\bigr)\E\bigl(\widehat
{D}_0(0,0)\bigr) (n\lambda)^{- 3/2},
\end{eqnarray*}
where we have used the fact that $\{\hat{D}_k(0,0) =
\hat{h}^k(\hat{r}(0,0))(1) - \hat{h}^k(\hat{l}(0,0))(1): k \geq0\}$
is a martingale (see Proposition \ref{propMartingale}).
From the earlier discussions, it also follows that $\E(\widehat
{D}_0(0,0))\leq
2\E(G) = 2(1-p)p^{-1}$ where $G$ is a geometric random variable. Thus,
$\limsup_{n\to\infty}\sqrt{n}\P(L(0,0)\leq n\tau, \# C(0,0) >
(n\lambda)^{3/2})
=0 $ as $\tau\to0$,
which completes the proof.
\end{pf}

\begin{pf*}{Proof of Theorem \ref{BEarea}}
We first recall the result Lemma 6.1 of Resnick \cite{R07},
page 174 which states that for nonnegative Radon measures $\mu, \mu
_n, n \geq1$, on
$[0,\infty)^d \setminus\{ \mathbf{0}\}$ we have $\mu_n\stackrel
{v}{\to} \mu$ if and only if $\mu_n  ([0,x_1] \times\cdots
\times[0, x_d] )^c \to
\mu ([0,x_1] \times\cdots\times[0, x_d] )^c$ for all
$x_1, \ldots, x_d \geq0$ with $(x_1, \ldots, x_d) \neq\mathbf{0}$.
This result implies that
Lemma \ref{lemmaLambdaBetaPositiveNonTruncate} together
with Corollary \ref{corLambdaPositiveNonTruncate}
and Theorem \ref{clusterheight} prove (\ref{eqnHack}).

Fix $\tau>0$, $\lambda> 0$. For $\alpha< 2/3$, $\delta> 0$ and
for all large $n$, we have
$\P(L(0,0)>n\tau,\#C(0,0)>(n\lambda)^{1/\alpha}) \leq
\P(L(0,0)>n\tau,\#C(0,0)>(n\delta)^{3/2})$. Fix any $\varepsilon> 0$
and choose
$\delta= \delta(\varepsilon) > 0$ so that $\frac{1}{\gamma_0\sqrt
{\pi\tau}}
\P ( \sqrt{2}\int_{0}^{\infty} W^{+,\tau}(t)\,dt > \wb{\delta
}  )
< \varepsilon$, where $\wb{\delta} = \delta^{3/2}(\gamma_0p)^{-1}$.
From Lemma \ref{lemmaLambdaBetaPositiveNonTruncate}, we have
\[
\limsup_{n \to\infty}\sqrt{n}\P\bigl(L(0,0)>n\tau,\#C(0,0)>(n\lambda
)^{1/\alpha}\bigr) < \varepsilon. %
\]

On the other hand,
from the properties of $W^{+}$ and $W^{\tau}$, it follows that
$\P(\int_0^\infty W^{+,\tau}(t)\,dt > 0) = 1$ for $\tau> 0$. Now for
$\alpha> 2/3$ and
$\delta> 0$ we have $\P(L(0,0)>n\tau,\#C(0,0) > (n\lambda
)^{1/\alpha}) \geq
\P(L(0,0)>n\tau,\#C(0,0) > (n\delta)^{3/2})$ for all large $n$.
Again from Lemma \ref{lemmaLambdaBetaPositiveNonTruncate}, we have
\begin{eqnarray*}
&& \liminf_{n \to\infty}\sqrt{n}\P\bigl(L(0,0)>n\tau,\#C(0,0)>(n
\lambda )^{1/\alpha}\bigr)
\\
&&\qquad \geq \frac{1}{\gamma_0\sqrt{\pi\tau}} \P \biggl( \sqrt{2}\int_{0}^{\infty}
W^{+,\tau}(t)\,dt > \wb{\delta } \biggr).
\end{eqnarray*}
Since
\begin{eqnarray*}
&& \limsup_{n \to\infty}\sqrt{n} \P\bigl(L(0,0)>n\tau,\#C(0,0)>(n
\lambda)^{1/\alpha}\bigr)
\\
&&\qquad  \leq\lim_{n \to\infty}\sqrt{n}\P\bigl(L(0,0)>n\tau\bigr) =
\frac
{1}{\gamma_0\sqrt{\pi\tau}},
\end{eqnarray*}
letting $\delta\to0$,
we have $\lim_{n \to\infty}\sqrt{n}\P(L(0,0)>n\tau,\#
C(0,0)>(n\lambda)^{1/\alpha}) =
\frac{1}{\gamma_0\sqrt{\pi\tau}}$ for $\alpha> 2/3$.
This completes the proof of (\ref{eqnTrivialCasesHack}).

The argument for $(L(0,0), (D_{\max}(0,0))^{1/2})$ being similar is
omitted. \end{pf*}

\section*{Acknowledgements}
Kumarjit Saha is grateful to the Indian Statistical Institute for
a fellowship to pursue his Ph.D. The authors also thank the referee for
comments which led to a significant improvement of the paper.




%

\printaddresses
\end{document}